\documentclass[11pt,a4paper,reqno]{amsart}
\usepackage[T1]{fontenc}
\usepackage{lmodern}
\usepackage{amsmath,amssymb,amsthm,mathtools,mathrsfs}
\usepackage{microtype,enumitem,hyperref}
\usepackage[margin=1in]{geometry}
\hypersetup{colorlinks=true,linkcolor=blue!55!black,citecolor=blue!55!black,urlcolor=blue!55!black}
\usepackage{xcolor}
\numberwithin{equation}{section}
\newtheorem{theorem}{Theorem}[section]
\newtheorem{proposition}[theorem]{Proposition}
\newtheorem{lemma}[theorem]{Lemma}
\newtheorem{corollary}[theorem]{Corollary}

\newcommand{\R}{\mathbb R}
\newcommand{\C}{\mathbb C}
\newcommand{\Sph}{\mathbb S^2}
\newcommand{\ip}[2]{\langle #1,#2\rangle}
\newcommand{\E}{\mathbb E}
\newcommand{\cH}{\mathcal H}
\newcommand{\cR}{\mathcal R}
\newcommand{\dd}{\,\mathrm d}
\newcommand{\one}{\mathbf 1}
\newcommand{\tr}{\operatorname{Tr}}

\newcommand{\supp}{\operatorname{supp}}
\title[Sharp Fourier extension inequalities on $
\mathbb{S}^2
$]{SHARP FOURIER EXTENSION INEQUALITIES ON $
\mathbb{S}^2
$ FOR ALL $q>4$}
\author{Hanli Tang}
\address[Hanli Tang]{Laboratory of Mathematics and Complex Systems (Ministry of Education), School of Mathematical Sciences, Beijing Normal University, Beijing, 100875, China}
\email{hltang@bnu.edu.cn}

\author{Zewei Xu}
\address[Zewei Xu]{Laboratory of Mathematics and Complex Systems (Ministry of Education), School of Mathematical Sciences, Beijing Normal University, Beijing, 100875, China}
\email{zwxu@mail.bnu.edu.cn}

\thanks{The first author was partly supported by the National Natural Science Foundation of China (Grant No.~12471100).}
\date{}
\begin{document}
\begin{abstract}
For the Fourier extension operator on $\mathbb S^2$, we determine the optimal $L^2(\mathbb S^2)\to L^q(\mathbb R^3)$ constant for every $q>4$, and characterize the nonzero complex extremizers. The result extends the sharp quartic theorem of Foschi \cite{Foschi} and the even-exponent results of Oliveira e Silva and Quilodr\'an \cite{OSQ} to exponents for which the usual convolution identities are unavailable. For every finite $q>4$ and also at $q=\infty$, the nonzero extremizers are precisely scalar multiples of characters $e^{ia\cdot\omega}$, for some $a\in\R^3$.
\end{abstract}
\maketitle

\section{Introduction}
The Fourier restriction problem asks whether the Fourier transform of an $L^p(\mathbb R^d)$ function can be meaningfully restricted to a curved hypersurface when $
1\le p\le 2
$. Its dual formulation is the Fourier extension problem. For the unit sphere $\mathbb S^{d-1}\subset\mathbb R^d$, the theorem of Tomas \cite{Tomas}, together with Stein's endpoint result (see \cite[Chapter~VIII]{Stein}), gives
\[
 \|\widehat{f\sigma_{d-1}}\|_{L^q(\R^d)}
 \leq C_{d,q}\|f\|_{L^2(\mathbb S^{d-1})},
 \qquad q\geq\frac{2(d+1)}{d-1}.
\]
The lower bound on $q$ is sharp in the $L^2$ extension scale. 

For a fixed admissible exponent, it is natural to ask what the least possible constant in the inequality is, whether it is attained, and which functions attain equality. We address these questions for the $L^2\to L^q$ Fourier extension inequality on $\mathbb S^2\subset\R^3$.

Let $\sigma$ denote the usual surface measure on $\Sph$. The Fourier extension of $f$ is given by
\begin{equation}\label{eq:Tdef}
 \widehat{f\sigma}(x)
 =\int_{\Sph}e^{-ix\cdot\omega}f(\omega)\dd\sigma(\omega).
\end{equation}
Motivated by the rotational symmetry of the sphere, we consider constant functions as candidates for extremizers. Building on this, if we define $\phi(x):=(4\pi)^{-1}\widehat{\sigma}(x)$, then
\begin{equation}\label{eq:phi}
 \phi(x)=\frac{\sin|x|}{|x|},\quad \phi(0)=1,\qquad
 A_q:=\int_{\R^3}|\phi(x)|^q\dd x
 =4\pi\int_0^\infty\left|\frac{\sin r}{r}\right|^q r^2\dd r.
\end{equation}
Taking $f\equiv1$, for which $\widehat{f\sigma}=4\pi\phi$ and $\|f\|_{L^2(\sigma)}=\sqrt{4\pi}$, we obtain the following lower bound for the best constant $C_q$ in the $L^2(\sigma)\to L^q(\R^3)$ extension inequality:
\[
 C_q\geq\sqrt{4\pi}\,A_q^{1/q},\qquad 4\leq q<\infty.
\]
Here and throughout the paper, a character on $\Sph$ means a function of the form $e^{ia\cdot\omega}$ for some $a\in\R^3$. In the context of spherical extension inequalities, multiplication by such a character corresponds to translation in $\R^3$: if $f(\omega)=ce^{ia\cdot\omega}$, then
$\widehat{f\sigma}(x)=4\pi c\phi(x-a)$ and $\|f\|_{L^2(\sigma)}=\sqrt{4\pi}|c|$.
Thus we ask whether the best constant is $C_q=\sqrt{4\pi}\,A_q^{1/q}$ and whether the family $\{ce^{ia\cdot\omega}:c\in\C\setminus\{0\},\ a\in\R^3\}$ contains all nonzero extremizers.

At the Stein--Tomas endpoint $q=4$ on $\Sph$, the sharp constant and the classification of extremizers were obtained through the work of Christ and Shao \cite{ChristShaoExist,ChristShaoExt} and Foschi \cite{Foschi}. Christ and Shao first established the existence of extremizers and compactness of nonnegative extremizing sequences, and later proved smoothness of critical points and precompactness of complex extremizing sequences modulo characters. They also showed that every complex extremizer is a nonzero scalar multiple of a nonnegative extremizer times a character, thereby reducing the classification to the nonnegative case.

Foschi \cite{Foschi} settled this remaining question. More precisely, Foschi proved  that the nonnegative extremizers are exactly the nonzero constant functions, and the optimal constant can be computed explicitly. Together with Christ and Shao's results \cite{ChristShaoExist,ChristShaoExt}, this gives the full classification of complex extremizers at $q=4$. 

A key feature of Foschi's proof is that the exponent is the even integer
$4=2\cdot2$. Writing
$f_\star(\omega):=\overline{f(-\omega)}$, Plancherel's theorem gives
\[
 \|\widehat{f\sigma}\|_{L^4(\R^3)}^4
 =(2\pi)^3\|f\sigma*f_\star\sigma\|_{L^2(\R^3)}^2.
\]
Thus the quartic norm is converted into the $L^2$ norm of a twofold
convolution. For the constant function, the relevant convolution density is
\[
 (\sigma*\sigma)(x)=\frac{2\pi}{|x|}\,
 \one_{\{0<|x|<2\}}(x).
\]
Since this convolution density is unbounded near the origin, a
direct application of Cauchy--Schwarz does not yield the sharp quartic
estimate. Foschi overcame this difficulty by exploiting a geometric
identity for four points on the sphere, together with the permutation
symmetry of the associated quadrilinear form, thereby handling the
singularity without losing sharpness.

We now record the result below.

\begin{theorem}[Foschi \cite{Foschi}]\label{thm:foschi}
For every $f\in L^2(\Sph,\sigma)$,
\begin{equation}\label{eq:foschi}
 \|\widehat{f\sigma}\|_{L^4(\R^3)}\leq 2\pi\|f\|_{L^2(\sigma)},\qquad A_4=\pi^2.
\end{equation}
The constant is sharp. Moreover, a nonzero function
$f\in L^2(\Sph,\sigma)$ is an extremizer for the inequality
\eqref{eq:foschi} if and only if
\begin{equation}\label{eq:foschi-extremizers}
 f(\omega)=ce^{ia\cdot\omega},\qquad
 c\in\C\setminus\{0\},\quad a\in\R^3.
\end{equation}
\end{theorem}

Building on Foschi's sharp result at exponent $q=4$ on $\Sph$,
Carneiro and Oliveira e Silva \cite{CarneiroOS} extended the sharp
quartic inequality to ambient dimensions $4\leq d\leq7$. Oliveira e
Silva and Quilodr\'an \cite{OSQ} then established the sharp
inequalities for every even exponent $q\geq6$ in dimensions
$d\in\{3,4,5,6,7\}$.

A major advantage of working with even exponents is that the
Plancherel identity continues to yield a multiple-convolution
formula. Indeed, for every even exponent $q=2m$, where $m\geq2$ is
an integer, we have
\begin{equation}\label{eq:even-convolution-intro}
 \|\widehat{f\sigma}\|_{L^{2m}(\R^3)}^{2m}
 =(2\pi)^3\big\|(f\sigma)^{*m}\big\|_{L^2(\R^3)}^2.
\end{equation}
Here $*$ denotes convolution of measures, and the identity holds
first for suitable functions and then by approximation. This
convolution structure provides a common starting point for the
sharp inequalities at all even exponents. In addition to the multiple-convolution identity, their arguments also exploit symmetrization, spectral information arising from the symmetries of the problem, and probabilistic properties of spherical convolutions. They proved that constant functions are the unique
real-valued maximizers and also characterized the complex-valued
maximizers. Specializing their results to $\Sph$ gives the
following theorem.

\begin{theorem}[Oliveira e Silva and Quilodr\'an \cite{OSQ}]\label{thm:osq-even}
Let $q\geq6$ be an even integer. Then, for every $f\in L^2(\Sph,\sigma)$,
\[
 \|\widehat{f\sigma}\|_{L^q(\R^3)}
 \leq\sqrt{4\pi}\,A_q^{1/q}\|f\|_{L^2(\sigma)}.
\]
The constant is sharp, and equality for a nonzero $f$ holds if and only if
$f(\omega)=ce^{ia\cdot\omega}$ for some $c\in\C\setminus\{0\}$ and
$a\in\R^3$.
\end{theorem}

Besides these even-exponent results, Oliveira e Silva and Quilodr\'an \cite{OSQ} also obtained an interesting conditional propagation result. More precisely, Theorem~1.5 of the same paper gives such a conclusion for sufficiently large $q$, provided that $d>2$ and constant functions maximize at exponent $q>q_*(d)$; under this assumption, every real-valued, continuously differentiable maximizer at exponent $q+2$ is constant. We recall its precise form on $\mathbb S^2$ in Section~\ref{sec:prelim}. Our reduction method developed in this paper is closely related to their conditional propagation result.

Beyond the sharp results at even exponents discussed above, Fanelli, Vega and Visciglia \cite[Theorem~0.1]{FVV} proved the existence of extremizers for the $L^2(\Sph,\sigma)\to L^q(\R^3)$ Fourier extension inequality for every finite real exponent $q>4$, including non-even exponents. They also proved that normalized maximizing sequences are precompact in $L^2(\sigma)$ after suitable modulations and passage to a subsequence. Thus, in the range considered here, existence is known. What remains is to determine the form of the maximizers.

A related but distinct extension problem is obtained by measuring the radial and angular variables separately in the target norm. Vega established the corresponding radial--angular mixed-norm inequality in his 1988 doctoral thesis \cite{Vega}, and Córdoba \cite{Cordoba} later gave an alternative proof:
\[
 \|\widehat{f\sigma_{d-1}}\|_{L^q_{\mathrm{rad}}L^2_{\mathrm{ang}}(\R^d)}
 \leq \mathbf C_{d,q}\|f\|_{L^2(\mathbb S^{d-1})},
 \qquad d\geq2,\quad q>\frac{2d}{d-1}.
\]
Carneiro, Oliveira e Silva and Sousa \cite{COS} subsequently studied its sharp form. Using spherical-harmonic decomposition, they reduced the sharp-constant problem to comparing weighted Bessel moments and proved that constant functions are the unique extremizers for all admissible even exponents, at $q=\infty$, and for all sufficiently large real exponents. Ciccone and Sousa \cite{CicconeSousa} later extended the range of exponents, covering the entire Stein--Tomas range in dimensions $d=2,3$. The second author of the present paper determined the sharp mixed-norm constant and characterized all nonzero extremizers for every dimension $d\geq2$ and every admissible real exponent $q>2d/(d-1)$ in \cite{XuMixed}. The same work also established a global quadratic stability inequality with the optimal coefficient for all these exponents, as well as for $q=\infty$. The mixed norm is more directly diagonalized by spherical harmonics than the $L^q(\R^3)$ norm considered here, so these results do not determine the sharp constant for the present problem.

Alongside the global sharp-constant and extremizer results discussed above, related work has also addressed local maximizers and stability questions for restriction inequalities. Gon\c{c}alves and Negro \cite{GoncalvesNegro} considered local maximizers for restriction inequalities associated with the sphere and other quadratic surfaces. Carneiro, Negro and Oliveira e Silva \cite{CarneiroNegroOS} proved sharp quartic inequalities and extremizer results under suitable small weighted perturbations. There are also recent results in different settings: Carneiro, Negro and Oliveira e Silva \cite{CarneiroNegroOSArbitrary} obtained sharp quartic inequalities in arbitrary dimensions with an additional term involving the spherical mean, while Wang \cite{WangCircle} determined the sharp endpoint constant on the circle. For more results concerning sharp Fourier extension inequalities, we refer the interested reader to \cite{FrankLiebSabin,FlockStovall,ShaoWang} and the references cited therein.

The results above leave a question in the non-even range. On $\mathbb S^2$, the sharp constant and all extremizers are known at $q=4$ and at every even integer $q\geq6$. Interpolation between the known even exponents gives an upper bound, but in general it does not give the candidate value $\sqrt{4\pi}\,A_q^{1/q}$. We prove that this candidate is in fact sharp for every $q>4$.

\begin{theorem}\label{thm:main}
For every finite $q>4$ and every $f\in L^2(\Sph,\sigma)$,
\begin{equation}\label{eq:main}
 \|\widehat{f\sigma}\|_{L^q(\R^3)}\leq C_q\|f\|_{L^2(\sigma)}.
\end{equation}
The constant is optimal and is given by
\begin{equation}\label{eq:area}
 C_q=\sqrt{4\pi}\,A_q^{1/q}.
\end{equation}
Moreover, a nonzero function $f\in L^2(\Sph,\sigma)$ is an extremizer
for the inequality \eqref{eq:main} if and only if
\begin{equation}\label{eq:chars}
 f(\omega)=ce^{ia\cdot\omega},\qquad
 c\in\C\setminus\{0\},\quad a\in\R^3.
\end{equation}
At $q=\infty$, the sharp inequality is
\[
 \|\widehat{f\sigma}\|_{L^\infty(\R^3)}
 \leq\sqrt{4\pi}\,\|f\|_{L^2(\sigma)}.
\]
Moreover, a nonzero function $f\in L^2(\Sph,\sigma)$ is an extremizer
for this $L^\infty$ inequality if and only if it has the form
\eqref{eq:chars}.
\end{theorem}

Let us give a brief overview of the main ideas in the proof of our theorem. The proof has two steps for finite $q$. We first reduce the problem for all $q>4$ to the interval $4<q\leq6$. The conditional propagation theorem of Oliveira e Silva and Quilodr\'an recalled in Section~\ref{sec:prelim} can be viewed as a precursor to the exponent-propagation method developed here. The key quantitative ingredient in this reduction is the following inequality:
\begingroup
\setlength{\abovedisplayskip}{7pt plus 1pt minus 2pt}
\setlength{\belowdisplayskip}{7pt plus 1pt minus 2pt}
\setlength{\abovedisplayshortskip}{7pt plus 1pt minus 2pt}
\setlength{\belowdisplayshortskip}{7pt plus 1pt minus 2pt}
\[
 \frac{A_{p+2}}{A_p}\geq\frac{p+1}{p+6},\qquad p\geq4,
\]
\endgroup
where equality holds only at $p=4$. The strict inequality for $p>4$ is what drives the propagation. Indeed, if $f$ is a normalized maximizer whose extension is
real-valued at exponent $p+2$ and is not a character, then
its Euclidean orbit supplies five independent tangent directions. Differentiating the Euler--Lagrange equation along these directions and using the trace estimate forces the $p$-th moment of $Tf$ to be strictly larger than $A_p$, contradicting the sharp estimate at $p$. Thus both sharpness and the characterization of extremizers propagate from $p$ to $p+2$.

In the second step, we show that every maximizer is a character for $4<q\leq6$. The existence theorem of Fanelli, Vega and Visciglia \cite[Theorem~0.1]{FVV} guarantees a maximizer, so we can argue by contradiction and suppose that there exists a maximizer that is not a character. We then derive a necessary condition from the Euler--Lagrange equation and its second variation. Thus, the rest of the proof is devoted to showing that this necessary condition cannot hold. Since the estimates depend on the exponent, we divide the interval into three parts and employ different estimation techniques, eventually obtaining a contradiction in all cases.

Our use of the second variation to derive a necessary
condition and rule out noncharacter maximizers is inspired
by the work of Frank and Lieb
\cite[Sections~3 and~5]{FrankLiebHeisenberg}
on sharp Hardy--Littlewood--Sobolev inequalities on the
Heisenberg group. We complement this variational argument
with the moment and spherical-harmonic estimates needed
to obtain the contradiction for the spherical Fourier
extension problem.

It is also worth noting that the combination of a sharp quartic estimate and a second-variation argument may be useful for non-endpoint extension inequalities on spheres in other dimensions, possibly for non-even exponents. We do not pursue this question here.

Section~\ref{sec:prelim} collects the basic properties of the extension operator, the spherical-harmonic decomposition, and the real-output reduction. We also recall the existence theorem of Fanelli, Vega and Visciglia and the conditional propagation theorem of Oliveira e Silva and Quilodr\'an, and derive the Euler--Lagrange equation for a maximizer. In Section~\ref{sec:reduction}, we prove the continuous moment inequality and combine it with the Euclidean orbit argument to reduce the problem to $4<q\leq6$. Section~\ref{sec:baseproof} then rules out noncharacter maximizers on this interval. The appendices contain the exact convolution moments and the spherical-harmonic estimates used in the proof.

\section{Preliminaries}\label{sec:prelim}
For the proofs, it is convenient to use the normalized surface measure
$\mu=\sigma/(4\pi)$ and the operator
\[
 Tf(x)=\int_{\Sph}e^{-ix\cdot\omega}f(\omega)\dd\mu(\omega)
 =\frac{1}{4\pi}\widehat{f\sigma}(x).
\]
Then $\|f\|_{L^2(\sigma)}=\sqrt{4\pi}\|f\|_{L^2(\mu)}$,
so the finite-exponent inequality in Theorem~\ref{thm:main} is equivalent to
$\|Tf\|_{L^q(\R^3)}\leq A_q^{1/q}\|f\|_{L^2(\mu)}$.
We first record two elementary properties of $T$. Since $\mu$ has
mass one, $\|Tf\|_\infty\leq\|f\|_2$. The Fourier transform of the
compactly supported finite measure $f\mu$ is real analytic. Its decay at
infinity requires a separate observation: for smooth $g$ on $\Sph$,
stationary phase on the sphere gives $Tg\in C_0(\R^3)$. Given
$f\in L^2(\mu)$, choose smooth $g_j\to f$ in $L^2(\mu)$. Then
\[
 \|Tf-Tg_j\|_\infty\leq\|f-g_j\|_{L^1(\mu)}
 \leq\|f-g_j\|_2\to0.
\]
Thus $Tf\in C_0(\R^3)$ for every $f\in L^2(\mu)$.
The uniqueness theorem for Fourier transforms of finite measures also
shows that $Tf\not\equiv0$ when $f\ne0$. In particular, the absolute
maximum of $Tf$ is attained and is positive for nonzero $f$.

We write $P_k$ for the orthogonal projection onto spherical harmonics of degree $k$, and $\mathcal H_k=P_kL^2(\mu)$, so that $\Delta_{\Sph}|_{\mathcal H_k}=-k(k+1)$. Radial weights on $\R^3$ lead to rotation-invariant quadratic forms and therefore do not mix different harmonic degrees. More precisely, the Funk--Hecke formula assigns to a zonal kernel $F(\omega\cdot\eta)$ the eigenvalue
\begin{equation}\label{eq:funk}
 \frac12\int_{-1}^1 F(t)\,\mathsf P_k(t)\dd t
\end{equation}
on $\mathcal H_k$, where $\mathsf P_k$ is the Legendre polynomial.

We recall the real-output reduction in \cite[Proposition~2.1]{OSQ}; see also \cite[Proposition~6.7]{BrocchiOSQ}. Let
\begin{equation}\label{eq:H}
 Jf(\omega)=\overline{f(-\omega)},\qquad
 \cH=\{f\in L^2(\mu):Jf=f\}.
\end{equation}
This is a real Hilbert space and $Tf$ is real-valued for $f\in\cH$. Every $f$ decomposes as $f=f_0+if_1$ with $f_0,f_1\in\cH$ and $\|f\|_2^2=\|f_0\|_2^2+\|f_1\|_2^2$. For $q>2$,
\begin{equation}\label{eq:realreduction}
 \|Tf\|_q^2=\|(Tf_0)^2+(Tf_1)^2\|_{q/2}
 \leq\|Tf_0\|_q^2+\|Tf_1\|_q^2.
\end{equation}
Thus the operator norm can be computed using $\cH$. The equality case is equally useful: strict convexity of $L^{q/2}$, analyticity, and uniqueness of Fourier transforms show that a nonzero complex extremizer differs by a constant phase from an extremizer in $\cH$. We will perform the variational argument in this real Hilbert space and return to complex functions at the end.

At the quartic endpoint, Theorem~\ref{thm:foschi} already identifies the sharp constant and all equality cases. For $q>4$, our argument starts instead from a maximizer and studies its Euler--Lagrange equation and second variation. We therefore need an existence result, together with compactness of maximizing sequences after removing the modulation symmetry. We use the following existence and compactness theorem.

\begin{theorem}[Fanelli, Vega and Visciglia \cite{FVV}]\label{thm:exist}
For every real $q>4$, the operator $T:L^2(\mu)\to L^q(\R^3)$ attains its norm. Its maximizing sequences are precompact in $L^2(\mu)$ after suitable modulations and extraction of subsequences.
\end{theorem}
For a unit real-output maximizer $f$ at an exponent $q>4$, the Euler--Lagrange equation reads
\begin{equation}\label{eq:EL}
 G_q(f)=Nf,\qquad G_q(f)=T^*(|Tf|^{q-2}Tf),\qquad N=\|Tf\|_q^q.
\end{equation}
Here $T^*:L^{q/(q-1)}(\R^3)\to L^2(\mu)$ denotes the adjoint of $T$.

For comparison with the work of Oliveira e Silva and Quilodr\'an, it is
convenient to introduce the functional
\[
 \Phi_q(f):=\frac{\|Tf\|_{L^q(\R^3)}^q}{\|f\|_{L^2(\mu)}^q},
 \qquad 0\ne f\in L^2(\mu).
\]
Their functional $\Phi_{3,q}$ is $(4\pi)^{q/2}\Phi_q$ in our
normalization, so the two functionals have the same maximizers. We recall
their conditional propagation result, specialized to dimension three.

\begin{theorem}[Oliveira e Silva and Quilodr\'an \cite{OSQ}]\label{thm:osq-bootstrap}
There exists a finite threshold $q_*(3)\geq4$ such that, whenever
$q>q_*(3)$ and constant functions maximize $\Phi_q$, every
real-valued $C^1$ maximizer of $\Phi_{q+2}$ is constant.
\end{theorem}

The finiteness of $q_*(3)$ follows from Theorem~1.5 of \cite{OSQ}.
The theorem is conditional: it assumes that constants already
maximize at exponent $q$, and it requires a $C^1$ maximizer at
exponent $q+2$. In Section~\ref{sec:reduction}, we establish a
moment comparison in a larger range and prove the regularity needed
for our orbit argument directly from the Euler--Lagrange equation.

\section{Reduction to the continuous interval $4<q\leq6$}\label{sec:reduction}
Theorem~\ref{thm:osq-bootstrap} gives a conditional propagation result
from $q$ to $q+2$ for $q>q_*(3)$, and assumes that the maximizer at
the higher exponent is $C^1$. We will use the same general principle
at every real exponent $p\geq4$. The scalar comparison behind the
result of Oliveira e Silva and Quilodr\'an is, in our normalization,
\[
 \frac{A_{p+2}}{A_p}>\frac{p+1}{p+6}.
\]
Their Proposition~7.8 establishes this comparison for $p>8$.
We first prove it for every real $p>4$ (with equality at $p=4$), and
then combine it with a trace estimate along the Euclidean orbit of a
maximizer. The proof of the moment comparison does not distinguish
integer from noninteger exponents.

\subsection*{A continuous moment-ratio inequality}
The propagation argument will require a comparison between two moments whose exponents differ by two. More precisely, the condition we need is
\begin{equation}\label{eq:ratio}
 \frac{A_{p+2}}{A_p}\geq\frac{p+1}{p+6},\qquad p\geq4.
\end{equation}
To prove this comparison, we integrate the radial Helmholtz equation for $\phi$ against a suitable power of $\phi$. The resulting identity involves the auxiliary function $Z$ appearing below. Its sign has to be coordinated with the size of $|\phi|$: the positive region should be precisely the region where $|\phi|$ dominates its later values. We first establish this sign pattern.

\begin{lemma}[Single sign change]\label{lem:sign}
With $\phi(r)=\sin r/r$ denoting the radial profile in \eqref{eq:phi},
put $Z(r)=\phi(r)^2+5\phi'(r)^2-1$. There is a unique
$R\in(0,\pi)$ such that
\begin{equation}\label{eq:Zsign}
 Z(r)>0\quad(0<r<R),\qquad Z(r)<0\quad(r>R),
 \qquad \phi(R)>\frac9{35}>\frac14\geq\sup_{r\geq\pi}|\phi(r)|.
\end{equation}
\end{lemma}
\begin{proof}
The radial Helmholtz equation is $\phi''+2\phi'/r+\phi=0$. Define
$H(r)=(2r^2-5)\sin r+5r\cos r$. Then
\begin{equation}\label{eq:Zderivative}
 Z'(r)=-\frac{4\phi'(r)}{r^3}H(r),\qquad H'(r)=r\sin r(2r\cot r-1).
\end{equation}
The function $r\cot r$ decreases strictly from $1$ to $-\infty$ on $(0,\pi)$, since
$(r\cot r)'=(\sin r\cos r-r)/\sin^2r<0$.
Thus $H'$ changes sign exactly once. Because $H(0)=0$ and $H(\pi)=-5\pi$, $H$ has one zero and changes from positive to negative. Also $\phi'<0$ on $(0,\pi)$, so $Z$ first increases and then decreases. As $Z(0)=0$ and $Z(\pi)=5/\pi^2-1<0$, it has exactly one zero $R$ in the first lobe.

For the height comparison, take $a=\sqrt6$. Alternating Taylor bounds give
$9/35\leq\phi(a)\leq3/10$ and $0<-\phi'(a)/a\leq529/3080<11/64$; hence
$Z(a)<9/100+30(11/64)^2-1<0$ and $R<a$. Monotonicity gives $\phi(R)>\phi(a)\geq9/35$.
For $r\geq\pi$, direct estimates yield
$\phi^2+5(\phi')^2\leq r^{-2}+5(r^2+1)r^{-4}\leq59/81<1$.
Finally, for $\pi\leq r\leq4$ we have $|\phi(r)|\leq(r-\pi)/r<1/4$, while for $r\geq4$ the bound $|\phi(r)|\leq1/r\leq1/4$ suffices.
\end{proof}
The sign change in Lemma~\ref{lem:sign} occurs before $|\phi|$ can exceed its first-lobe value again. Hence increasing the moment exponent favors exactly the region where $Z$ is positive. The following argument makes this observation quantitative.
\begin{theorem}\label{thm:ratio}
Inequality \eqref{eq:ratio} holds for every real $p\geq4$, with equality only at $p=4$.
\end{theorem}
\begin{proof}
We first derive the identity behind the comparison. Since $\phi(r)=\sin r/r$,
\[
 (r^2\phi'(r))'=-r^2\phi(r).
\]
Multiply this equation by $\phi|\phi|^p$ and integrate on $(\varepsilon,L)$.
The derivative of $\phi|\phi|^p$ is $(p+1)|\phi|^p\phi'$ for every
real $p>1$, including at zeros of $\phi$. An integration by parts therefore gives
\[
 \begin{split}
 (p+1)\int_\varepsilon^L |\phi|^p(\phi')^2r^2\dd r
 &=\int_\varepsilon^L|\phi|^{p+2}r^2\dd r\\
 &\quad+\bigl[r^2\phi'(r)\phi(r)|\phi(r)|^p\bigr]_{\varepsilon}^{L}.
 \end{split}
\]
The boundary term vanishes as $\varepsilon\to0^+$ because
$\phi'(r)=O(r)$ there; at infinity it is $O(L^{-p})$, since both
$\phi$ and $\phi'$ are $O(L^{-1})$. We obtain
\begin{equation}\label{eq:helmholtz}
 (p+1)\int_0^\infty |\phi|^p(\phi')^2r^2\dd r=\frac{A_{p+2}}{4\pi}.
\end{equation}
In particular, multiplying $Z=\phi^2+5(\phi')^2-1$ by
$|\phi|^pr^2$ and using this identity gives
\begin{equation}\label{eq:ratioZ}
 \frac{(p+6)A_{p+2}-(p+1)A_p}{4\pi}
 =(p+1)\int_0^\infty|\phi|^pZ(r)r^2\dd r.
\end{equation}
It remains to determine the sign of the last integral. Taking $p=4$
in \eqref{eq:ratioZ} and using the exact moments
$A_4=\pi^2$ and $A_6=\pi^2/2$, we obtain
\[
 \frac{10A_6-5A_4}{4\pi}
 =5\int_0^\infty|\phi(r)|^4Z(r)r^2\dd r=0.
\]
Thus the integral in question vanishes at $p=4$. For $p>4$, we subtract
this zero integral multiplied by $\phi(R)^{p-4}$, where $R$ is the
unique zero from Lemma~\ref{lem:sign}. This gives
\[
 \int_0^\infty|\phi(r)|^pZ(r)r^2\dd r
 =\int_0^\infty|\phi(r)|^4
 \bigl(|\phi(r)|^{p-4}-\phi(R)^{p-4}\bigr)Z(r)r^2\dd r.
\]
The signs of the two factors agree on either side of $R$.
If $0<r<R$, then $r<\pi$ and the strict decrease of $\phi$ on
$(0,\pi)$ gives
\[
 \phi(r)>\phi(R)>0.
\]
Since $p-4>0$, it follows that
$|\phi(r)|^{p-4}-\phi(R)^{p-4}>0$. At the same time,
$Z(r)>0$ by Lemma~\ref{lem:sign}. Hence the product of these two
factors is positive on $(0,R)$.

If $R<r<\pi$, the same monotonicity gives
$|\phi(r)|=\phi(r)<\phi(R)$. If $r\geq\pi$, Lemma~\ref{lem:sign}
gives
\[
 |\phi(r)|\leq\frac14<\frac9{35}<\phi(R).
\]
Consequently, $|\phi(r)|^{p-4}-\phi(R)^{p-4}<0$ for every $r>R$.
Since $Z(r)<0$ there, the product is again positive. At $r=R$ both
factors vanish. We have therefore shown that the integrand on the
right-hand side is nonnegative everywhere and strictly positive on
$(0,R)$; in particular, its integral is strictly positive. Equation
\eqref{eq:ratioZ} now proves the strict inequality for $p>4$, while
$p=4$ is the asserted equality case.
\end{proof}

\subsection*{Regularity and the Euclidean orbit}
The trace argument requires differentiating rotation orbits in $L^2$, which is justified once the maximizer belongs to $H^1(\Sph)$. Related regularity arguments appear in \cite{ChristShaoExt,OSQSmooth}; we give the proof for the real exponents needed here.
\begin{lemma}[Orbit regularity]\label{lem:H1}
Let $2<q<\infty$ and suppose $T:L^2(\mu)\to L^q(\R^3)$ is bounded.
If a nonzero $f\in L^2(\mu)$ satisfies $G_q(f)=\lambda f$, where
$G_q$ is defined in \eqref{eq:EL}, then $\lambda>0$ and $f\in H^1(\Sph)$.
\end{lemma}
\begin{proof}
Pairing $G_q(f)=\lambda f$ with $f$ gives
\[
 \lambda\|f\|_2^2=\langle G_q(f),f\rangle=\|Tf\|_q^q.
\]
Since $f\mu$ is a finite measure and the Fourier transform is injective
on finite measures, $Tf$ cannot vanish identically when $f\ne0$.
Thus $\|Tf\|_q^q>0$ and $\lambda>0$.

We view the relevant Banach spaces as real spaces when differentiating.
For $\mathcal N_q(u)=|u|^{q-2}u$, the real derivative is
\[
 D\mathcal N_q(u)v
 =|u|^{q-2}v+(q-2)|u|^{q-4}u\operatorname{Re}(\bar u v),
 \quad |D\mathcal N_q(u)v|\leq(q-1)|u|^{q-2}|v|,
\]
with derivative zero at $u=0$. The corresponding Nemytskii map is
continuously Fr\'echet differentiable from $L^q$ to $L^{q/(q-1)}$
for $q>2$ (see \cite{AppellZabrejko}). In particular,
$G_q=T^*\mathcal N_qT$ is $C^1$ as a real map on $L^2(\mu)$.
All adjoints and pairings in the compactness argument below are
understood over the underlying real Hilbert spaces.

To differentiate the identity along a rotation orbit, we first need an
$L^2$ derivative of that orbit. The compactness of the linearized
Euler--Lagrange operator provides it. Fix a ball $B_R\subset\R^3$ and
write
\[
 T_{B_R}:L^2(\mu)\to L^2(B_R),\qquad
 T_{B_R}h=(Th)|_{B_R}.
\]
This operator is Hilbert--Schmidt, since
\[
 \|T_{B_R}\|_{\mathrm{HS}}^2
 =\int_{B_R}\int_{\Sph}|e^{-ix\cdot\omega}|^2
       \dd\mu(\omega)\dd x
 =|B_R|.
\]
Moreover, $\|Tf\|_\infty\leq\|f\|_2$, and hence
\[
 \|D\mathcal N_q(Tf(x))\|_{\mathcal L_{\R}(\C,\C)}
 \leq(q-1)\|Tf\|_\infty^{q-2}.
\]
Thus pointwise application of $D\mathcal N_q(Tf(x))$ defines a bounded
real-linear multiplication operator on $L^2(B_R)$. Set
\[
 DG_{q,R}(f)
 :=T_{B_R}^*
   \bigl(D\mathcal N_q(Tf)|_{B_R}\bigr)T_{B_R},
\]
where the middle factor denotes this multiplication operator.
Since $T_{B_R}$ is Hilbert--Schmidt, $DG_{q,R}(f)$ is compact.

We next show that these compact operators converge to $DG_q(f)$ in
operator norm. For $h,k\in L^2(\mu)$, the pointwise derivative bound
and H\"older's inequality give
\[
 \begin{split}
 |\langle (DG_q(f)-DG_{q,R}(f))h,k\rangle|
 &\leq(q-1)\int_{\R^3\setminus B_R}
       |Tf|^{q-2}|Th|\,|Tk|\,\dd x\\
 &\leq(q-1)\|Tf\|_{L^q(\R^3\setminus B_R)}^{q-2}
       \|Th\|_q\|Tk\|_q\\
 &\leq(q-1)\|T\|_{2\to q}^2
       \|Tf\|_{L^q(\R^3\setminus B_R)}^{q-2}
       \|h\|_2\|k\|_2.
 \end{split}
\]
The second line uses H\"older with exponents $q/(q-2),q,q$.
Since $Tf\in L^q(\R^3)$, the tail norm tends to zero as $R\to\infty$.
Consequently,
\[
 \|DG_q(f)-DG_{q,R}(f)\|_{2\to2}\to0.
\]
As an operator-norm limit of compact operators, $DG_q(f)$ is compact.

Let $U_t$ be a one-parameter rotation group on $L^2(\mu)$ and put
$h_t=U_tf-f$. The equation is rotation-equivariant, so $U_tf$
satisfies the same equation as $f$. Since $h_t\to0$ in $L^2$,
Fr\'echet differentiability gives
\begin{equation}\label{eq:H1-diff}
 \lambda h_t=DG_q(f)h_t+r_t,\qquad
 \|r_t\|_2=o(\|h_t\|_2).
\end{equation}
Let $P_N$ be the orthogonal projection onto the sum of all spherical
harmonic spaces of degrees at most $N$. Then $P_N\to I$ strongly.
Since $DG_q(f)$ is compact, this strong convergence implies
\[
 \|(I-P_N)DG_q(f)\|_{2\to2}\to0.
\]
Choose $N$ so that, with $P=P_N$,
\[
 \|(I-P)DG_q(f)\|_{2\to2}\leq\frac{\lambda}{4}.
\]
Rotations preserve each spherical-harmonic degree, so $P$ commutes
with $U_t$. The vector $Pf$ belongs to a finite-dimensional space
of smooth functions; hence its rotation orbit is differentiable and
\[
 \|Ph_t\|_2=\|U_tPf-Pf\|_2\leq C|t|.
\]
Applying $I-P$ to \eqref{eq:H1-diff}, we obtain
\[
 \lambda\|(I-P)h_t\|_2
 \leq\frac{\lambda}{4}\|h_t\|_2+o(\|h_t\|_2).
\]
For sufficiently small $t$, the remainder is at most
$\frac{\lambda}{4}\|h_t\|_2$. It follows that
\[
 \|(I-P)h_t\|_2\leq\frac12\|h_t\|_2.
\]
Combining this with the estimate for $Ph_t$ gives
\[
 \|h_t\|_2
 \leq\|Ph_t\|_2+\|(I-P)h_t\|_2
 \leq C|t|+\frac12\|h_t\|_2,
\]
and therefore
\[
 \|U_tf-f\|_2\leq 2C|t|.
\]
Thus the rotation difference quotients are bounded in $L^2(\mu)$.

For $1\leq i<j\leq3$, let
\[
 X_{ij}(\omega)=\omega_i e_j-\omega_j e_i,\qquad
 D_{ij}=X_{ij}\cdot\nabla_{\Sph}
       =\omega_i\partial_{\omega_j}-\omega_j\partial_{\omega_i}
\]
be the corresponding rotation vector field and its differential
operator. Choose the one-parameter rotations $R_{ij}(t)$ whose
velocity field is $X_{ij}$, and set
$U_tg(\omega)=g(R_{ij}(t)\omega)$. For any sequence $t_n\to0$,
the difference quotients
\[
 \frac{U_{t_n}f-f}{t_n}
\]
are bounded in $L^2$, so a subsequence converges weakly to some
$g\in L^2(\mu)$. For every smooth test function $\psi$, unitarity
of the rotation action gives
\[
 \left\langle\frac{U_{t_n}f-f}{t_n},\psi\right\rangle
 =\left\langle f,\frac{U_{-t_n}\psi-\psi}{t_n}\right\rangle
 \to-\langle f,D_{ij}\psi\rangle.
\]
Hence
\[
 \langle g,\psi\rangle=-\langle f,D_{ij}\psi\rangle,
\]
which says precisely that $D_{ij}f=g$ in the distributional sense.
Repeating this argument for the three coordinate rotations shows
that $D_{ij}f\in L^2(\mu)$ for every $i<j$; the subsequences may be
chosen separately for the three generators.

Finally, the vector fields $X_{ij}$ span the tangent plane at every
point of $\Sph$. More precisely, for every smooth function $v$,
\[
 \sum_{1\leq i<j\leq3}|D_{ij}v(\omega)|^2
 =|\omega|^2|\nabla_{\Sph}v(\omega)|^2
   -|\omega\cdot\nabla_{\Sph}v(\omega)|^2
 =|\nabla_{\Sph}v(\omega)|^2.
\]
Since the fields $X_{ij}$ span the tangent plane, on each
coordinate patch the coordinate derivatives are smooth linear
combinations of the $D_{ij}$. Thus the distributional coordinate
derivatives of $f$ belong locally to $L^2$, which gives
$f\in H^1(\Sph)$. By density, the displayed identity extends to
these weak derivatives as well.
\end{proof}
The Sobolev regularity just obtained allows us to differentiate the full Euclidean symmetry action, including rotations, at a maximizer.
For $(a,R)\in E^+(3)=\R^3\rtimes SO(3)$, define
\[
 (\rho(a,R)f)(\omega)=e^{ia\cdot\omega}f(R^{-1}\omega),\qquad
 T(\rho(a,R)f)(x)=Tf(R^{-1}(x-a)).
\]
The action preserves the real-output space $\cH$ defined in
\eqref{eq:H} and both norms. If $f\in\cH$ is a unit maximizer, let
$u=Tf$, $N=\|u\|_q^q$ and $K_f=T^*|u|^{q-2}T$. Its generator in the direction $X=(v,A)\in\R^3\rtimes\mathfrak{so}(3)$ is
\begin{equation}\label{eq:generator}
 L_Xf=i(v\cdot\omega)f-(A\omega)\cdot\nabla_{\Sph}f\in L^2(\mu).
\end{equation}
Differentiating the Euler--Lagrange equation along the orbit gives
\begin{equation}\label{eq:tangent}
 \ip{L_Xf}{f}=0,\qquad K_f(L_Xf)=\frac{N}{q-1}L_Xf.
\end{equation}
These identities show that every nonzero direction tangent to the Euclidean orbit is an eigenvector of $K_f$ with the same eigenvalue. For the trace estimate, however, the number of independent such directions matters. The next lemma gives the required lower bound for the dimension of the orbit of a noncharacter.

\begin{lemma}[Orbit dimension]\label{lem:orbit}
If nonzero $f\in\cH$ is not a scalar multiple of a character, then
its Euclidean orbit has dimension at least five.
\end{lemma}
\begin{proof}
Let $u=Tf\in C_0(\R^3)$, and let
\[
G_f=\{g\in E^+(3):\rho(g)f=f\}
\]
be the stabilizer of $f$. The representation $\rho$ is strongly
continuous on $L^2(\mu)$, so $G_f$ is closed. By the closed subgroup
theorem, it is a Lie subgroup of $E^+(3)$; see \cite{Lee}.

We first prove that $G_f$ is compact. Choose $x_0\in\R^3$ with
$u(x_0)\neq0$. The covariance identity
\[
T(\rho(g)f)(x)=u(g^{-1}x)
\]
shows that $|u(gx_0)|=|u(x_0)|$ for every $g\in G_f$. Since
$u\in C_0(\R^3)$, the superlevel set
\[
K=\{x\in\R^3:|u(x)|\geq |u(x_0)|\}
\]
is compact, and hence $G_fx_0\subset K$. Write $g=(a,R)$, so
$gx_0=a+Rx_0$ and therefore
\[
a=gx_0-Rx_0.
\]
The first term on the right ranges over the compact set $K$, while
the second ranges over the compact sphere of radius $|x_0|$.
Thus the translation components of elements of $G_f$ are bounded.
Their rotational components lie in $SO(3)$. Since $G_f$ is closed
in $\R^3\rtimes SO(3)$, it follows that $G_f$ is compact.

We can now average the orbit of $x_0$ with respect to normalized Haar
measure on $G_f$. Set
\[
b=\int_{G_f}g x_0\,dg.
\]
For every $h\in G_f$, left invariance of Haar measure gives
\[
hb=\int_{G_f}(hg)x_0\,dg
   =\int_{G_f}g x_0\,dg=b.
\]
Thus $b$ is a common fixed point of the affine action of $G_f$.
Let $\tau_b=(b,I_3)$ denote translation by $b$, and set
\[
\widetilde f=\rho(\tau_{-b})f=e^{-ib\cdot\omega}f.
\]
The stabilizer of $\widetilde f$ is
\[
G_{\widetilde f}=\tau_{-b}G_f\tau_b\subset SO(3),
\]
because every element of $G_f$ fixes $b$.

Suppose that the Lie algebra of $G_{\widetilde f}$ were all of
$\mathfrak{so}(3)$. Its identity component would then be $SO(3)$,
and $\widetilde f$ would be invariant under every rotation. The only
rotation-invariant functions in $L^2(\Sph)$ are constants. Hence
$\widetilde f$ would be constant, and
$f=e^{ib\cdot\omega}\widetilde f$ would be a scalar multiple of a
character, contrary to the hypothesis. Therefore the stabilizer Lie
algebra is a proper Lie subalgebra of $\mathfrak{so}(3)$.

Under the identification $\mathfrak{so}(3)\simeq\R^3$, the Lie
bracket corresponds to the cross product. If a subspace contained
two linearly independent vectors $v,w$, closure under the bracket
would force it to contain $v\times w$, which is nonzero and
orthogonal to both $v$ and $w$. Thus no two-dimensional subspace
of $\mathfrak{so}(3)$ is a Lie subalgebra. It follows that
\[
\dim G_f=\dim G_{\widetilde f}\leq1.
\]
Since $\dim E^+(3)=6$, the orbit, viewed as the homogeneous space
$E^+(3)/G_f$, has dimension
\[
\dim(E^+(3)\cdot f)=6-\dim G_f\geq5.
\]
\end{proof}

For the trace argument, we need the orbit directions to be genuine
$L^2$ tangent vectors, not merely infinitesimal directions in the
Lie-group quotient. We now use the additional fact that $f$ is a
maximizer. It satisfies the Euler--Lagrange equation
\eqref{eq:EL}, so Lemma~\ref{lem:H1} gives $f\in H^1(\Sph)$.
Consequently the orbit map
\[
\Phi_f:E^+(3)\to L^2(\mu),\qquad
\Phi_f(g)=\rho(g)f,
\]
is differentiable at the identity, and its differential is
\[
d\Phi_f(e)X=L_Xf,
\]
with $L_Xf$ as in \eqref{eq:generator}. The kernel of this
differential is precisely the Lie algebra of $G_f$. Indeed, if
$X$ belongs to that Lie algebra, the corresponding one-parameter
subgroup fixes $f$, so $L_Xf=0$. Conversely, if $L_Xf=0$, then
the differentiable curve $t\mapsto\rho(\exp(tX))f$ has zero
derivative at every $t$ and is therefore constant; hence
$\exp(tX)\in G_f$ for all $t$. The rank of $d\Phi_f(e)$ is thus
\[
6-\dim G_f\geq5.
\]
We may choose five linearly independent orbit tangents. By
\eqref{eq:tangent}, they all belong to $f^\perp$ and are eigenvectors
of $K_f$ with the same eigenvalue $N/(q-1)$. After orthonormalizing
their span, we obtain five orthonormal orbit tangents with the same
eigenvalue. Together with the radial direction $f$, these are the
six orthonormal vectors used in the trace estimate below.
\begin{proposition}[Strict trace obstruction]\label{prop:trace}
Let $q>4$, and suppose $T:L^2(\mu)\to L^{q-2}(\R^3)$ is bounded. If a unit real-output maximizer $f$ at $q$ is not a character, then
\begin{equation}\label{eq:trace}
 \int_{\R^3}|Tf|^{q-2}\dd x>\frac{q+4}{q-1}\int_{\R^3}|Tf|^q\dd x.
\end{equation}
\end{proposition}
\begin{proof}
Because $T:L^2(\mu)\to L^{q-2}(\R^3)$ is bounded,
$|u|^{q-2}\in L^1(\R^3)$. The weighted extension operator
$|u|^{(q-2)/2}T:\cH\to L^2(\R^3;\R)$ is Hilbert--Schmidt.
Indeed, for a real orthonormal basis $(e_j)$ of $\cH$, Parseval
and the unit norm of the evaluation functional $f\mapsto Tf(x)$ give
\[
 \tr_{\cH}K_f
 =\sum_j\int_{\R^3}|u(x)|^{q-2}|Te_j(x)|^2\dd x
 =\int_{\R^3}|u(x)|^{q-2}\dd x.
\]
In particular $K_f$ is positive and trace class.
The unit vector $f$ contributes $N$. By Lemma~\ref{lem:orbit} and the differentiability argument immediately following it, five orthonormal orbit tangents contribute $5N/(q-1)$ through \eqref{eq:tangent}. An additional unit vector $h$ orthogonal to those six vectors contributes strictly positively: $u$ is nonzero on a ball and the analytic function $Th$ cannot vanish throughout that ball. Thus
$\tr K_f>N+5N/(q-1)=(q+4)N/(q-1)$.
\end{proof}
The trace estimate becomes decisive when the exponent is raised by two: its left-hand side is then a moment of order $p$, which the known sharp $L^p$ estimate bounds from above.
\begin{corollary}[Propagation by two]\label{cor:bootstrap}
Suppose $p\geq4$ and $\|T\|_{2\to p}=A_p^{1/p}$, with only characters as nonzero extremizers. Then the same sharpness and rigidity hold at $p+2$.
\end{corollary}
\begin{proof}
Existence at $q=p+2$ is Theorem~\ref{thm:exist}. Let $f\in\cH$ be a normalized maximizer, $N=\|Tf\|_q^q$. By the sharp $p$-estimate,
$\tr K_f=\int|Tf|^p\leq A_p$.
If $f$ is not a character, Proposition~\ref{prop:trace} implies
\[
 A_p>\frac{p+6}{p+1}N\geq\frac{p+6}{p+1}A_{p+2}\geq A_p,
\]
where the last inequality is Theorem~\ref{thm:ratio}. This contradiction proves the result; \eqref{eq:realreduction} supplies the complex equality cases.
\end{proof}
Unlike Theorem~\ref{thm:osq-bootstrap}, the present propagation result
applies to every real $p\geq4$. This extension uses the moment comparison
throughout that range, and Lemma~\ref{lem:H1} supplies the $L^2$
differentiability of the Euclidean orbit without assuming that the
maximizer is $C^1$.
Since the propagation argument holds for every real $p\geq4$, it can be iterated with increments of two. Given any $q>4$, we may subtract an even integer and reach an exponent in $(4,6]$. The following proposition records this reduction.
\begin{proposition}[Base-interval reduction]\label{prop:base}
To prove Theorem~\ref{thm:main} for all finite real $q>4$, it suffices to prove it for $4<q\leq6$.
\end{proposition}
\begin{proof}
For every finite $q>4$, choose an integer $n\geq0$ such that $q-2n\in(4,6]$. If $n=0$, the assertion is part of the base-interval result. If $n\geq1$, apply Corollary~\ref{cor:bootstrap} successively at $q-2n,q-2n+2,\ldots,q-2$. No interpolation between integer exponents is required.
\end{proof}

\section{Proof on the base interval $4<q\leq6$}\label{sec:baseproof}
Proposition~\ref{prop:base} leaves the initial interval $4<q\leq6$. Suppose that the sharp inequality fails there. The compactness theorem produces a maximizer which is not a character. After a modulation, its extension has a positive maximum at the origin, and we may separate the constant spherical harmonic from the remaining part of the data.

The second variation in that remainder leads to a necessary inequality
$\cR_p(r,v^2)\geq1$ (Proposition~\ref{prop:necessary}). The estimates below show that no admissible parameters can satisfy it. We treat the three subranges from the upper endpoint down. On $5\leq q\leq6$, the exact even-moment data give a direct interpolation argument. Below $q=5$, a fractional-moment estimate is needed to control the constant mode. As $q$ approaches $4$, this separate control is no longer sufficient, and the degree-two contribution must be retained in both the quartic deficit and the spectral estimate. Thus the order of the argument follows the increasing strength of the estimates required; The subdivision reflects the estimates available in each range.

\subsection{The second-variation necessary condition}\label{subsec:necessity}
To turn this setup into a usable obstruction, we combine the orthogonal decomposition with first- and second-order information at the peak. The first-order condition removes the degree-one component, while a normalized perturbation in the remaining direction supplies the second-variation constraint. The next proposition states the resulting decomposition and scalar inequality.
\begin{proposition}[Necessary condition for a noncharacter maximizer]\label{prop:necessary}
Let $4<q\leq6$, and suppose that $f\in\cH$ (as defined in \eqref{eq:H}),
$\|f\|_2=1$, is a maximizer not proportional to a character. Modulate and change its sign so that $Tf$ takes its positive absolute maximum at the origin. Then
\begin{equation}\label{eq:peakdecomp}
 f=c+\rho H,\qquad c,\rho>0,\quad c^2+\rho^2=1,\quad
 \|H\|_2=1,\quad P_0H=P_1H=0.
\end{equation}
Write $W=TH$, $N=\|Tf\|_q^q$, $p=q-2$, and $v=\rho/c$, and define $r,b_p\geq0$ by
\begin{equation}\label{eq:rbdef}
 r^p=\frac{\int_{\R^3}|\phi|^pW^2\dd x}{N},\qquad
 b_p^p=(1+v^2)^{(p-2)/2}\left(v^2+\frac1{p+1}\right).
\end{equation}
Then
\begin{equation}\label{eq:necessary}
 1\leq\cR_p(r,v^2):=
 \frac{p(p+1)rv^2}{(1+v^2)^{p/2}-1}
 \int_0^1(1-t)\bigl((1-t)r+tb_p\bigr)^{p-1}\dd t.
\end{equation}
\end{proposition}
\begin{proof}
Choose the origin at a point where $|Tf|$ is maximal. Since $Tf\in C_0(\R^3)$,
its absolute maximum is attained. Modulate the input to bring a point of
maximum to the origin, and change its sign so that $Tf(0)=c>0$.
Because $\mu$ is a probability measure and $f\in\cH$,
\[
 c=\int_{\Sph}f\dd\mu=\langle f,\one\rangle_{L^2(\mu)},
 \qquad 0<c\leq\|f\|_2=1.
\]
If $c=1$, equality in Cauchy--Schwarz would make $f$ constant, contrary
to our hypothesis. Thus $0<c<1$. Orthogonally decompose
$f=c+\rho H$, where $\rho=\sqrt{1-c^2}>0$, $\|H\|_2=1$ and
$P_0H=0$. The first derivatives at the maximum vanish:
\[
 0=\partial_{x_j}Tf(0)=-i\int_{\Sph}\omega_j f(\omega)\dd\mu(\omega)
 =-i\rho\int_{\Sph}\omega_jH(\omega)\dd\mu(\omega).
\]
The coordinate functions span the first spherical harmonics; this
proves $P_1H=0$ and establishes \eqref{eq:peakdecomp}.

The character input has norm one and extension norm $A_q^{1/q}$,
so extremality gives $N\geq A_q$. On the other hand,
$\|Tf\|_\infty=c$ by construction, whence
\[
 N=\int|Tf|^{q-4}|Tf|^4
 \leq c^{q-4}\|Tf\|_4^4\leq c^{q-4}A_4.
\]
We have obtained
\begin{equation}\label{eq:peakbound}
 A_q\leq N\leq c^{q-4}A_4.
\end{equation}

We vary the nonconstant component in the direction $H$, keeping the $L^2$ norm equal to one. The curve
\[
 f_\tau=\frac{f+\tau H}{(1+2\rho\tau+\tau^2)^{1/2}}
\]
consists of unit vectors for small real $\tau$. Let $u=Tf$ and
$W=TH$, and write $B(\tau)=\int|u+\tau W|^q$. As $u,W$ are real-valued,
\[
 B(0)=N,\quad B'(0)=q\int|u|^{q-2}uW,\quad
 B''(0)=q(q-1)\int|u|^{q-2}W^2.
\]
Differentiating $B(\tau)(1+2\rho\tau+\tau^2)^{-q/2}$ at zero,
maximality first gives
\[
 \int|u|^{q-2}uW=\rho N.
\]
The second derivative at the maximum is nonpositive. Using the
first-derivative identity to collect the terms from normalization gives
\[
 \left.\frac1q\frac{\dd^2}{\dd\tau^2}\|Tf_\tau\|_q^q\right|_{\tau=0}
 =(q-1)\int|u|^{q-2}W^2\dd x
   -\bigl(1+(q-2)\rho^2\bigr)N\leq0.
\]
 Using $u=c(\phi+vW)$ and
$c^2=(1+v^2)^{-1}$, one obtains
\[
 \int|\phi+vW|^pW^2\dd x
 \leq N(1+v^2)^{(p-2)/2}
 \left(v^2+\frac1{p+1}\right)=Nb_p^p.
\]

To relate this spectral information to the constant component, compare the Euler--Lagrange equation at $f$ with the corresponding expression at the constant function.
Define
\[
 F(t)=\int_{\R^3}|\phi+tvW|^p(\phi+tvW)\phi\dd x
 \qquad(0\leq t\leq1).
\]
Because the radial function $|\phi|^p\phi$ is orthogonal, under
integration on $\R^3$, to the extension of every harmonic of positive
degree, we have
\[
 F(0)=A_q,\qquad F'(0)=(p+1)v\int|\phi|^p\phi W=0.
\]
Pair the Euler--Lagrange equation $G_q(f)=Nf$ with $\one$.
As $\langle f,\one\rangle=c$ and $u=c(\phi+vW)$,
\[
 c^{p+1}F(1)=Nc,\qquad F(1)=N(1+v^2)^{p/2}.
\]
The peak estimate $A_q\leq N$ therefore implies
\[
 \begin{split}
 N\bigl((1+v^2)^{p/2}-1\bigr)
 &\leq F(1)-F(0)\\
 &=\int_0^1(1-t)F''(t)\dd t.
 \end{split}
\]
For real $p\geq2$ we may differentiate the scalar function
$y\mapsto|y|^py$ twice, including at $y=0$, and find
\[
 F''(t)=p(p+1)v^2\int
 |\phi+tvW|^{p-2}(\phi+tvW)\phi W^2\dd x.
\]
H\"older in the positive measure $W^2\dd x$ gives
\[
 |F''(t)|\leq p(p+1)v^2
 \left(\int|\phi|^pW^2\right)^{1/p}
 \left(\int|\phi+tvW|^pW^2\right)^{(p-1)/p}.
\]
The first factor equals $N^{1/p}r$. Minkowski and the second-variation
estimate bound the second norm by
\[
 \begin{split}
 \left(\int|\phi+tvW|^pW^2\right)^{1/p}
 &\leq(1-t)\left(\int|\phi|^pW^2\right)^{1/p}
  +t\left(\int|\phi+vW|^pW^2\right)^{1/p}\\
 &\leq N^{1/p}\bigl((1-t)r+tb_p\bigr).
 \end{split}
\]
Substituting these estimates and canceling $N>0$ gives
\eqref{eq:necessary}. All integrands in the differentiations are
locally dominated by $C(|\phi|+v|W|)^q$, which belongs to $L^1$;
this is also why zeros of $\phi+tvW$ cause no difficulty.
\end{proof}
The remaining estimates concern only $r$ and $v$. Every noncharacter maximizer would have to satisfy \eqref{eq:necessary}; we shall show that the sharp quartic inequality and the harmonic spectrum prevent this.

For use below put
\begin{equation}\label{eq:eta}
 S=1+v^2,\qquad \eta_p=1-\frac{p}{(p+1)S},\qquad
 b_p=\sqrt S\,\eta_p^{1/p},\qquad x=\frac r{b_p}.
\end{equation}
The integral in \eqref{eq:necessary} can be evaluated directly;
for $0\leq x<1$,
\[
 \int_0^1(1-t)\bigl((1-t)x+t\bigr)^{p-1}\dd t
 =\frac{1-(p+1)x^p+px^{p+1}}{p(p+1)(1-x)^2}.
\]
Consequently, the variational expression becomes
\begin{equation}\label{eq:Fp}
 \cR_p(r,S-1)=\frac{r(S-1)b_p^{p-1}}{S^{p/2}-1}\,F_p(x),\qquad
 F_p(x)=\frac{1-(p+1)x^p+px^{p+1}}{(1-x)^2}.
\end{equation}
In particular,
\begin{equation}\label{eq:F234}
 F_2(x)=1+2x,\quad F_3(x)=1+2x+3x^2,\quad
 F_4(x)=1+2x+3x^2+4x^3.
\end{equation}
For fixed $0<x<1$, $F_p(x)$ increases with $p$. This follows by
differentiating in $p$ and using $-\log x\geq1-x$. We will interpolate these functions between $p=2,3,4$ while keeping the variational parameters $r$ and $S$ fixed. The following lemma gives the needed uniform estimate.

\begin{lemma}[Two endpoint comparisons]\label{lem:interpolateR}
For the following scalar estimates, keep the spectral parameter $r$
fixed when writing
$K_p(S)=(1-S^{-1})/(1-S^{-p/2})$ (with $K_p(1)=2/p$) and
$\mathcal A_j(p,S)=r\sqrt S\,\eta_p^{1-1/p}F_j(r/b_p)$.
Then
$\cR_p(r,S-1)=K_p(S)\mathcal A_p(p,S)$.
For $p=n+t$, $n\in\{2,3\}$, $0\leq t\leq1$, suppose $1\leq S\leq S_*$, $r/b_p<1$, and
$\mathcal A_j(n,S)\leq M_j$ for $j=n,n+1$. Then
\begin{equation}\label{eq:Rinterp}
 \cR_{n+t}(r,S-1)\leq h_n(t)
 \bigl(K_n(S_*)M_n\bigr)^{1-t}
 \bigl(K_{n+1}(S_*)M_{n+1}\bigr)^t,
\end{equation}
where
\begin{equation}\label{eq:hn}
 h_n(t)=\frac{(n+t)(n+t+1)}{n(n+1)((n+2)/n)^t},
 \qquad h_2(t)<\frac{16}{15},\quad h_3(t)<\frac{21}{20}.
\end{equation}
Moreover, $\cR_p(r,S-1)$ is increasing in $r$.
\end{lemma}
\begin{proof}
We first record the dependence of $K_p$ on $S$ and $p$. For $S>1$,
differentiation gives
\[
 \frac{\partial}{\partial S}\log K_p(S)
 =\frac1S\left(\frac1{S-1}-\frac{p/2}{S^{p/2}-1}\right)\geq0,
\]
where the inequality follows from
$S^{p/2}-1\geq(p/2)(S-1)$ for $p\geq2$. Moreover,
\[
 \frac{\partial^2}{\partial p^2}\log K_p(S)
 =\frac{(\log S)^2 S^{p/2}}{4(S^{p/2}-1)^2}>0.
\]
These properties extend to $S=1$ by continuity, since
$K_p(1)=2/p$. Thus, whenever $p=n+t$ and $0\leq t\leq1$,
\begin{equation*}
 K_p(S)\leq K_n(S)^{1-t}K_{n+1}(S)^t
 \leq K_n(S_*)^{1-t}K_{n+1}(S_*)^t.
\end{equation*}

Next, the integral representation
\[
 \frac{F_p(x)}{p(p+1)}
 =\int_0^1(1-u)\bigl((1-u)x+u\bigr)^{p-1}\dd u
\]
allows H\"older's inequality to be applied with exponents
$1/(1-t)$ and $1/t$, with the endpoint cases interpreted directly.
For $p=n+t$ it yields
\begin{equation*}
 F_{n+t}(x)\leq
 \frac{(n+t)(n+t+1)}{n(n+1)((n+2)/n)^t}
 F_n(x)^{1-t}F_{n+1}(x)^t
 =h_n(t)F_n(x)^{1-t}F_{n+1}(x)^t.
\end{equation*}
This is the interpolation step; the remaining point is to compare
the endpoint factors at the same values of $r$ and $S$.

For $j\in\{2,3,4\}$, the polynomial formulas in \eqref{eq:F234}
give
\[
 \mathcal A_j(p,S)
 =\sum_{m=1}^j m r^m S^{1-m/2}\eta_p^{1-m/p}.
\]
Here $r$ and $S$ remain fixed. Since
$\eta_p=(1+p(1-S^{-1}))/(p+1)$, we have
\[
 \frac{\partial}{\partial p}\log\bigl(\eta_p^{1-m/p}\bigr)
 =-\int_{1-S^{-1}}^1\frac{1+mu}{(1+pu)^2}\dd u<0.
\]
(The integral is positive also at $S=1$.) Each coefficient in the
last sum is therefore nonincreasing with $p$, and hence
$\mathcal A_j(n+t,S)\leq\mathcal A_j(n,S)$.
Combining this fact with the H\"older estimate and the bound for
$K_p$ gives
\[
 \begin{split}
 \cR_{n+t}(r,S-1)
 &\leq h_n(t)K_{n+t}(S)
    \mathcal A_n(n+t,S)^{1-t}
    \mathcal A_{n+1}(n+t,S)^t\\
 &\leq h_n(t)
    \bigl(K_n(S_*)M_n\bigr)^{1-t}
    \bigl(K_{n+1}(S_*)M_{n+1}\bigr)^t,
 \end{split}
\]
which proves \eqref{eq:Rinterp}.

For completeness, the numerical bounds on $h_n$ can be checked
without optimization. The inequalities $\log2>2/3$ and
$\log(5/3)>1/2$, together with $e^y\geq1+y+y^2/2$, give
\[
 2^t\geq1+\frac{2t}{3}+\frac{2t^2}{9},\qquad
 \left(\frac53\right)^t\geq1+\frac t2+\frac{t^2}{8}
 \quad(0\leq t\leq1).
\]
Using
\[
 h_2(t)=\frac{1+5t/6+t^2/6}{2^t},\qquad
 h_3(t)=\frac{1+7t/12+t^2/12}{(5/3)^t},
\]
the claimed estimates reduce respectively to
\[
 \frac{19t^2-33t+18}{270}>0,\qquad
 \frac{23t^2-28t+24}{480}>0.
\]
Both quadratic polynomials are strictly positive (their
discriminants are negative). Finally, the integral expression in
\eqref{eq:necessary}, or its continuous extension at $S=1$, is
nondecreasing in $r\geq0$ for $p\geq2$. This proves the last assertion.
\end{proof}
The next lemma handles the range $1\leq S\leq2$ and $r\leq1/\sqrt5$, which will recur in the two lower-exponent cases.
\begin{lemma}[Small peaks]\label{lem:smallpeaks}
For $2\leq p\leq3$, $1\leq S\leq2$, and $0\leq r\leq1/\sqrt5$,
\begin{equation}\label{eq:smallpeaks}
 \cR_p(r,S-1)<1.
\end{equation}
At $S=1$, the expression is understood by its continuous extension,
as in Lemma~\ref{lem:interpolateR}.
\end{lemma}
\begin{proof}
Since $\cR_p(r,S-1)$ increases with $r$, it suffices to take
$R=1/\sqrt5$. Before applying Lemma~\ref{lem:interpolateR}, we check
that its condition $r/b_p<1$ is satisfied. Differentiating the formula
for $b_p$ in \eqref{eq:eta} gives
\[
 \frac{\partial}{\partial p}\log b_p
 =\int_{1-S^{-1}}^1\frac{u}{(1+pu)^2}\dd u\geq0.
\]
Consequently, for $2\leq p\leq3$,
\[
 b_p\geq b_2=\sqrt{S-\tfrac23}
 \geq\frac1{\sqrt3}>R.
\]

Write $w=b_2=\sqrt{S-2/3}$. As $1\leq S\leq2$,
\[
 \frac1{\sqrt3}\leq w\leq\frac2{\sqrt3}.
\]
Using the explicit polynomials in \eqref{eq:F234}, the two endpoint
quantities in Lemma~\ref{lem:interpolateR} are
\[
 \mathcal A_2(2,S)=Rw+2R^2,\qquad
 \mathcal A_3(2,S)=Rw+2R^2+\frac{3R^3}{w}.
\]
The first is increasing in $w$. The second is convex, since its
second derivative is $6R^3/w^3>0$, so its maximum occurs at an
endpoint. Moreover, the difference between its values at the upper
and lower endpoints is
\[
 \mathcal A_3(2,2)-\mathcal A_3(2,1)
 =\frac{R}{\sqrt3}\left(1-\frac{9R^2}{2}\right)
 =\frac{R}{10\sqrt3}>0.
\]
Thus both endpoint quantities are bounded by their respective values
at $S=2$:
\[
 \mathcal A_2(2,S)\leq M_2:=\frac25+\frac2{\sqrt{15}}
 <\frac{23}{25},\qquad
 \mathcal A_3(2,S)\leq M_3:=\frac25+\frac{29}{10\sqrt{15}}
 <\frac76.
\]

Now write $p=2+t$, $0\leq t\leq1$. Applying
\eqref{eq:Rinterp} with $n=2$ and $S_*=2$ yields
\[
 \cR_{2+t}(R,S-1)
 \leq h_2(t)\bigl(K_2(2)M_2\bigr)^{1-t}
                  \bigl(K_3(2)M_3\bigr)^t.
\]
Here $K_2(2)=1$ and
$K_3(2)=(4+\sqrt2)/7<11/14$. In particular,
\[
 K_2(2)M_2<\frac{23}{25},\qquad
 K_3(2)M_3<\frac{11}{14}\cdot\frac76
 =\frac{11}{12}<\frac{23}{25}.
\]
Since $h_2(t)<16/15$ by \eqref{eq:hn}, we conclude that
\[
 \cR_p(r,S-1)
 \leq\cR_p(R,S-1)
 <\frac{16}{15}\cdot\frac{23}{25}
 =\frac{368}{375}<1.
\]
This proves the lemma throughout the stated parameter range.
\end{proof}
For $5\leq q\leq6$, we will only obtain $S\leq8/3$ and $r\leq4/9$. The following comparison covers that range.
\begin{lemma}[Upper scalar comparison]\label{lem:upperpeaks}
For $3\leq p\leq4$, $1\leq S\leq8/3$, and $0\leq r\leq4/9$,
\begin{equation}\label{eq:upperpeaks}
 \cR_p(r,S-1)<1.
\end{equation}
\end{lemma}
\begin{proof}
By monotonicity in $r$ it suffices to take $R=4/9$. In the
following derivative check, $\eta_3=1-3/(4S)$ and
$x=R/(\sqrt S\eta_3^{1/3})$ are exactly the quantities already
defined in \eqref{eq:eta} at $p=3$, $r=R$.
After factoring out the positive expression
$R\sqrt S\eta_3^{2/3}/(2S(S-3/4))$, the derivatives of
$\mathcal A_3(3,S)$ and $\mathcal A_4(3,S)$ have the respective signs of
\[
 \begin{gathered}
 S+\tfrac14-3R^2\eta_3^{1/3}+x,\\
 S+\tfrac14-3R^2\eta_3^{1/3}
       -\frac{8R^3}{\sqrt S}+x(1-2x^2).
 \end{gathered}
\]
Both are positive: $\eta_3^{1/3}<9/10$, $x<1/\sqrt2$, and even
the second expression exceeds $5/4-(27/10)R^2-8R^3>0$.
Thus both endpoint quantities are largest at $S=8/3$, where
$8/9<\eta_3^{1/3}<9/10$ and
\[
 \mathcal A_3(3,8/3)<\frac{1361}{1215},\qquad
 \mathcal A_4(3,8/3)<\frac{1441}{1215}.
\]
As $K_3(8/3)<5/6$, $K_4(8/3)=8/11$ and
$(8/11)(1441/1215)<(5/6)(1361/1215)$, \eqref{eq:Rinterp} yields
\[
 \cR_p(r,S-1)<\frac{21}{20}\cdot\frac56\cdot\frac{1361}{1215}<1.
\]
\end{proof}

\subsection{Spectral diagonalization and endpoint data}\label{subsec:common-spectral}
We now turn to the weighted quadratic form that defines $r$:
\[
 H\longmapsto\int_{\R^3}|\phi(x)|^{q-2}|TH(x)|^2\dd x.
\]
Its associated positive self-adjoint operator is
$T^*|\phi|^{q-2}T$. Since $\phi$ is radial, this operator commutes
with rotations. By the Funk--Hecke formula, it acts as a scalar on
each spherical-harmonic space $\mathcal H_k$ introduced in
Section~\ref{sec:prelim}. We denote this scalar
by $\beta_{q,k}$. Equivalently, for any $Y_k\in\mathcal H_k$ with
$\|Y_k\|_2=1$,
\[
 \beta_{q,k}=\int_{\R^3}|\phi(x)|^{q-2}|TY_k(x)|^2\dd x.
\]
Thus $\beta_{q,k}$ depends only on $q$ and the harmonic degree $k$,
not on the choice of $Y_k$; in particular, $\beta_{q,0}=A_q$.
The exact convolution calculations in Appendix~\ref{app:convol}
give the following values:

\begin{equation}\label{eq:spectra}
\begin{aligned}
 A_4&=\pi^2,\qquad A_6=\frac{\pi^2}{2},\qquad A_8=\frac{\pi^2}{3},\\
 \beta_{4,k}&=\frac{\pi^2}{2k+1},\qquad
 \beta_{6,k}=\frac{3\pi^2/2}{(2k-1)(2k+1)(2k+3)}\quad(k\geq1).
\end{aligned}
\end{equation}

The scalar action on each $\mathcal H_k$ gives the orthogonal
expansion needed below. For $H=\sum_{k\geq0}P_kH$, we have
\begin{equation}\label{eq:spectral-parseval}
 \int_{\R^3}|\phi|^{q-2}|TH|^2\dd x
 =\sum_{k\geq0}\beta_{q,k}\|P_kH\|_2^2.
\end{equation}
One may first use finite sums and then pass to general $H$ by
boundedness of the quadratic extension operator. In particular, the degree-two contribution and the sum over degrees $k\geq3$ can be estimated separately.

For the first subrange, only the endpoint eigenvalues in
\eqref{eq:spectra} and this decomposition are needed. The fractional-moment estimates for $q<5$ are developed in Subsection~\ref{subsec:middle}.

\subsection{The range $5\leq q\leq6$: endpoint interpolation}\label{subsec:highrange}
We first treat the subrange nearest $q=6$, where the exact even-exponent
moments and spectra suffice. The aim is to show that a hypothetical
noncharacter maximizer has parameters $r$ and $S=1+v^2$
(defined in \eqref{eq:rbdef} and \eqref{eq:eta}) satisfying
the hypotheses of Lemma~\ref{lem:upperpeaks}. We bound the denominator
$N$ in the definition of $r$ from below, then its numerator by
interpolating the fourth- and sixth-order quadratic forms. The same
lower bound for $N$, combined with the peak estimate, controls $S$.

\par\medskip
\noindent\emph{Constant-mode moment.}\quad
The map $q\mapsto\log A_q$ is convex by H\"older's inequality.
For $5\leq q\leq6$, convexity places the graph above the line
through the points at $q=6$ and $q=8$, extended to the left:
\[
 \log A_q\geq\log A_6
      +\frac{q-6}{2}\bigl(\log A_8-\log A_6\bigr).
\]
Inserting $A_6/A_4=1/2$ and $A_8/A_4=1/3$ from
\eqref{eq:spectra}, we obtain
\begin{equation}\label{eq:easyA}
 \frac{A_q}{A_4}\geq\frac34(2/3)^{(q-4)/2}.
\end{equation}
This lower bound will be used both for the spectral ratio $r$
and for the peak coefficient $c$.

\par\medskip
\noindent\emph{Spectral interpolation.}\quad
For $H$ in \eqref{eq:peakdecomp}, $\|H\|_2=1$ and
$P_0H=P_1H=0$. Thus only spherical-harmonic degrees $k\geq2$
contribute. The exact values in \eqref{eq:spectra} imply
\[
 \frac{\beta_{4,k}}{A_4}=\frac1{2k+1}\leq\frac15,
 \qquad
 \frac{\beta_{6,k}}{\beta_{4,k}}
 =\frac{3}{2(2k-1)(2k+3)}\leq\frac1{14}.
\]
By the orthogonal expansion \eqref{eq:spectral-parseval} and
$\sum_{k\geq2}\|P_kH\|_2^2=1$, these yield the two endpoint bounds
\[
 \int_{\R^3}|\phi|^2|TH|^2\dd x\leq\frac{A_4}{5},
 \qquad
 \int_{\R^3}|\phi|^4|TH|^2\dd x\leq\frac{A_4}{70}.
\]
To interpolate between them, put $t=(q-4)/2\in[0,1]$.
H\"older's inequality, applied with the common factor $|TH|^2$,
gives
\[
 \begin{split}
 \int_{\R^3}|\phi|^{q-2}|TH|^2\dd x
 &=\int_{\R^3}
   (|\phi|^2|TH|^2)^{1-t}
   (|\phi|^4|TH|^2)^t\dd x\\
 &\leq\left(\int_{\R^3}|\phi|^2|TH|^2\dd x\right)^{1-t}
       \left(\int_{\R^3}|\phi|^4|TH|^2\dd x\right)^t\\
 &\leq\frac{A_4}{5}\,14^{-t}.
 \end{split}
\]
Consequently,
\begin{equation}\label{eq:easybeta}
 \frac1{A_4}\int_{\R^3}|\phi|^{q-2}|TH|^2\dd x
 \leq\frac15\,14^{-(q-4)/2}.
\end{equation}

\par\medskip
\noindent\emph{Bounds for $r$ and $S$.}\quad
We now use $p=q-2\in[3,4]$ as in Proposition~\ref{prop:necessary}.
A normalized maximizer has $N=\|Tf\|_q^q\geq A_q$.
By the definition of $r$ in \eqref{eq:rbdef},
\eqref{eq:easybeta}, and the lower bound \eqref{eq:easyA},
\[
 \begin{split}
 r^p
 &=\frac{\displaystyle\int_{\R^3}|\phi|^p|TH|^2\dd x}{N}\\
 &\leq\frac{(A_4/5)14^{-(p-2)/2}}{A_q}
 \leq\frac4{15}\left(\frac3{28}\right)^{(p-2)/2}.
 \end{split}
\]
Taking the $p$-th root and collecting the constant factors gives
\begin{equation}\label{eq:easy-r}
 r\leq\sqrt{\frac3{28}}\left(\frac{112}{45}\right)^{1/p}
 \leq\left(\frac4{525}\right)^{1/6}<\frac49.
\end{equation}
Here the middle inequality uses $p\geq3$ and $112/45>1$;
the last one follows, on taking sixth powers, from
$531441<537600$.

To control $S$, we instead use the other half of the peak estimate
\eqref{eq:peakbound}, namely $A_q\leq N\leq c^{q-4}A_4$.
Combining it with \eqref{eq:easyA} gives
\[
 c^{p-2}\geq\frac34\left(\frac23\right)^{(p-2)/2}.
\]
As $S=c^{-2}$ and $p-2\geq1$, we obtain
\begin{equation}\label{eq:easy-S}
 S\leq\frac32(4/3)^{2/(p-2)}\leq\frac83,
 \qquad v^2=S-1\leq\frac53.
\end{equation}
Thus the moment and spectral bounds provide, respectively, the
required restrictions on the peak parameter $S$ and the radius $r$.

Fix $q\in[5,6]$ and suppose that there exists a normalized
real-output maximizer $f\in\cH$ which is not a character.
By Proposition~\ref{prop:necessary}, we may modulate $f$ so that
its parameters satisfy the necessary condition \eqref{eq:necessary}.
In particular, $v>0$ and $S=1+v^2>1$. The bounds just proved give
\[
 3\leq p=q-2\leq4,\qquad
 1<S\leq\frac83,\qquad 0<r<\frac49.
\]
These are precisely the hypotheses of Lemma~\ref{lem:upperpeaks}.
That lemma therefore gives
\[
 \cR_p(r,v^2)=\cR_p(r,S-1)<1,
\]
whereas the second-variation necessary condition
\eqref{eq:necessary} asserts that the same quantity is at least one.
This contradiction rules out noncharacter real-output maximizers
throughout $5\leq q\leq6$.

Theorem~\ref{thm:exist} ensures that a maximizer exists at each such
$q$, and \eqref{eq:realreduction} allows the norm to be attained in
$\cH$. Every real-output maximizer must therefore be a character;
since characters have extension norm $A_q^{1/q}$, this is the optimal
constant. Finally, the equality case in \eqref{eq:realreduction}
shows that every nonzero complex maximizer is also a scalar multiple
of a character. This proves Theorem~\ref{thm:main} on $[5,6]$.

\subsection{The range $22/5\leq q\leq5$: a convexity bridge}\label{subsec:middle}
The endpoint interpolation above is sufficient for $p=q-2\geq3$,
but its spectral bound becomes too weak below $q=5$. We first
refine the degree-two eigenvalue using its first three moments;
the remaining harmonic degrees can still be handled by endpoint
interpolation. For the constant mode, a logarithmic moment provides
an additional fractional-exponent anchor. Together these estimates
will control the spectral radius and the peak parameter.

\subsubsection*{Fractional spectral and logarithmic moments}
As $q\to4^+$, the degree-two eigenvalue gives the largest contribution among the admissible harmonic degrees. To improve the estimate based on two endpoint moments, we use a two-node quadrature rule determined by three moments.

\begin{lemma}[Two-node moment comparison]\label{lem:gauss}
Let $0<X\leq1$ be a random variable whose law has more than two support points, and set
\[
 m_j=\E X^j\ (j=1,2,3),\qquad V=m_2-m_1^2,\quad
 U=\frac{m_3-m_1m_2}{V},\quad W=\frac{m_1m_3-m_2^2}{V}.
\]
Put $x_\pm=(U\pm\sqrt{U^2-4W})/2$ and $w=(x_+-m_1)/(x_+-x_-)$. Then $0<x_-<x_+<1$, $0<w<1$, and
\begin{equation}\label{eq:gauss}
 \E X^s\leq w x_-^s+(1-w)x_+^s\leq m_1^s,\qquad0\leq s\leq1.
\end{equation}
Both inequalities are strict for $0<s<1$.
\end{lemma}
\begin{proof}
Let $\nu$ be the law of $X$. Since $\supp\nu$ contains more than two
points, the inner product
\[
 (P,Q)_\nu=\int P(x)Q(x)\,\dd\nu(x)
\]
is positive definite on the polynomials of degree at most two.
Consequently there is a unique monic quadratic polynomial orthogonal
to $1$ and $x$. The two orthogonality relations are
\[
 m_2-Um_1+W=0,\qquad m_3-Um_2+Wm_1=0,
\]
and solving them gives
\[
 P_2(x)=x^2-Ux+W
\]
with precisely the coefficients in the statement.

The two-node Gaussian quadrature theorem
\cite[Theorems~1.45--1.46, §1.4.2]{Gautschi} now gives two distinct
nodes $x_-<x_+$ in the interior of the support interval and positive
weights $\lambda_-,\lambda_+$ such that
\[
 \E Q(X)=\lambda_-Q(x_-)+\lambda_+Q(x_+),
 \qquad \deg Q\leq3.
\]
The nodes are the zeros of $P_2$, hence they are the stated $x_\pm$.
Applying exactness to $Q=1$ and $Q(x)=x$ shows that
$\lambda_-+\lambda_+=1$ and
$\lambda_-x_-+\lambda_+x_+=m_1$. Thus
\[
 \lambda_-=\frac{x_+-m_1}{x_+-x_-}=w,\qquad
 \lambda_+=1-w.
\]
In particular, $0<x_-<x_+<1$ and $0<w<1$.

Fix $0<s<1$, and let $H_s$ be the cubic Hermite interpolant to
$x\mapsto x^s$ at $x_-$ and $x_+$, matching both the values and the
first derivatives. For each $x>0$, the Hermite interpolation remainder
formula (see \cite[Theorem~1.48 and the subsequent remainder formula
in §1.4.2]{Gautschi}) gives, unless $x$ is one of the nodes,
\[
 x^s-H_s(x)
 =\frac{s(s-1)(s-2)(s-3)}{4!}\,
   \xi_x^{\,s-4}(x-x_-)^2(x-x_+)^2
\]
for some $\xi_x>0$ between $x$ and the two nodes. This pointwise use
is legitimate even if $0$ belongs to the support interval: for each
fixed $x>0$, the function $t^s$ is $C^4$ on the compact interval
spanned by $x,x_-,x_+$. Since
\[
 s(s-1)(s-2)(s-3)<0,
\]
we have $H_s(x)\geq x^s$, with strict inequality away from the two
nodes. Exactness of the quadrature rule for the cubic polynomial $H_s$
therefore yields
\[
 \E X^s<\E H_s(X)
 =wH_s(x_-)+(1-w)H_s(x_+)
 =wx_-^s+(1-w)x_+^s.
\]
The strictness follows because the law of $X$ is not concentrated on
the two nodes. Finally, strict concavity of $x^s$ and the first-moment
identity give
\[
 wx_-^s+(1-w)x_+^s
 <\bigl(wx_-+(1-w)x_+\bigr)^s=m_1^s.
\]
The endpoint cases $s=0,1$ follow directly.
\end{proof}

Apply Lemma~\ref{lem:gauss} to the measure obtained by weighting physical space with $\phi^2|TY|^2$, where $Y$ is a unit degree-two harmonic. The associated random variable is $X=\phi^2$, and its first three moments can be calculated from the exact convolution formulas of Appendix~\ref{app:convol}.

For any unit $Y\in\mathcal H_2\cap\cH$, define the probability measure
\begin{equation}\label{eq:nuB}
 \dd\nu_B(x)=\frac{\phi(x)^2|TY(x)|^2}{\beta_{4,2}}\dd x,
 \qquad X(x)=\phi(x)^2.
\end{equation}
The convolution formula in Appendix~\ref{app:convol} gives
\begin{equation}\label{eq:three-moments}
 \E_{\nu_B}X=\frac1{14},\qquad
 \E_{\nu_B}X^2=\frac5{252},\qquad
 \E_{\nu_B}X^3=\frac{41}{4752}.
\end{equation}
Consequently the polynomial and weights in Lemma~\ref{lem:gauss} are
\begin{equation}\label{eq:PB}
 P_B(x)=x^2-\frac{1679}{3432}x+\frac{311}{20592},
 \quad w=\frac{x_+-1/14}{x_+-x_-}.
\end{equation}
The exact comparisons
\[
 \frac{331}{10000}<x_-<\frac1{30},\qquad
 \frac{57}{125}<x_+<\frac12,
 \qquad w>\frac{10}{11}
\]
follow by substituting these four rational points into $P_B$ and using
$57/125+10(331/10000)-11/14=9/7000>0$. Thus for $0\leq s\leq1$,
\begin{equation}\label{eq:Gbar}
 \E_{\nu_B}X^s\leq G(s):=wx_-^s+(1-w)x_+^s
 \leq\overline G(s):=\frac{10}{11}30^{-s}+\frac1{11}2^{-s}.
\end{equation}
We shall write $G$ for the exact two-node expression; it is log-convex with $G(0)=1$.
For $s>0$,
$x_-^s<x_+^s$; increasing the weight of the smaller node can only
reduce the weighted moment. Since $w>10/11$ and
$x_-<1/30$, $x_+<1/2$, replacing the nodes by their respective upper
bounds and replacing $w$ by $10/11$ gives the stated majorant.
 Simple fifth- and second-power comparisons give
\begin{equation}\label{eq:Gend}
 \overline G(1/5)<\frac{27}{50},\qquad
 \overline G(1/2)<\frac7{30}.
\end{equation}
The simpler bound $\overline G$ also controls all degrees $k\geq3$, allowing one estimate for the degree-two and higher-degree contributions.

For higher harmonic degrees, H\"older interpolation between \eqref{eq:spectra} yields, for $0\leq s\leq1$ and $k\geq3$,
\begin{equation}\label{eq:highmodes}
 \frac{\beta_{4+2s,k}}{A_4}\leq\frac1{2k+1}
 \left(\frac{3}{2(2k-1)(2k+3)}\right)^s
 \leq\frac17\,30^{-s}\leq\frac15\,\overline G(s).
\end{equation}
Combining degree two with $k\geq3$, we have the common estimate
\begin{equation}\label{eq:allspectral}
 \int|\phi|^{2+2s}|TH|^2\dd x\leq\frac{A_4}{5}\,\overline G(s)
 \quad(0\leq s\leq1,\ \|H\|_2=1,\ P_0H=P_1H=0).
\end{equation}

\subsubsection*{A logarithmic moment for the constant mode}
The denominator in the definition of $r$ is $N=\|Tf\|_q^q\geq A_q$. Near $q=4$ we need a quantitative lower bound for $A_q/A_4$. The logarithmic moment of $\phi^2$, under the fourth-power weight, provides such a bound by Jensen's inequality. Let
\[
 \dd\nu_4=A_4^{-1}\phi^4\dd x,\qquad X=\phi^2.
\]
The following logarithmic moment is the input for Jensen's inequality in
the near-endpoint range:
\begin{equation}\label{eq:entropy}
 \E_{\nu_4}X=\frac12,\quad\E_{\nu_4}X^2=\frac13,\quad
 \E_{\nu_4}(-\log X)=1+2\log2-2\gamma<\frac{37}{30}.
\end{equation}

\subsubsection*{A fractional moment at $q=22/5$}
The logarithmic moment in \eqref{eq:entropy}, together with the first two ordinary moments, allows a sharper lower bound for $A_{22/5}$. We obtain it by a Hermite minorant for $x^{1/5}$ in the span of $1,x,x^2,\log x$.
\begin{lemma}[Logarithmic Hermite minorant]\label{lem:logHermite}
For $0<s<1$ and $0<a<b$, there is a unique
$Q(x)=A+Bx+Cx^2+D\log x$ agreeing with $x^s$ in value and first derivative at $a$ and $b$. Its logarithmic coefficient is positive and
$Q(x)<x^s$ for all $x>0$ other than $a,b$.
\end{lemma}
\begin{proof}
The interpolation problem consists of four linear equations for the
four coefficients of $Q$. To see that their coefficient matrix is
nonsingular, consider a function
$A+Bx+Cx^2+D\log x$ with a double zero at each of $a,b$.
Four zeros counted with multiplicity give a zero of its third derivative
by repeated Rolle's theorem. But the third derivative is $2D/x^3$;
therefore $D=0$. A quadratic cannot have two distinct double zeros
unless it is identically zero. This establishes uniqueness and existence.

Now put $h(x)=x^s-Q(x)$. By construction $h(a)=h'(a)=h(b)=h'(b)=0$.
A repeated Rolle argument produces a point between $a$ and $b$
where $h'''=0$. Since
\[
 x^3h'''(x)=s(s-1)(s-2)x^s-2D,
\]
and $s(s-1)(s-2)>0$ for $0<s<1$, this forces $D>0$.
In particular $h(x)\to+\infty$ as $x\to0^+$, because the
term $-D\log x$ dominates the bounded terms and $x^s$.

For completeness, the same differentiation argument also rules out
an additional zero of $h$. Indeed
\[
 \frac{\dd}{\dd x}\bigl(x^3h'''(x)\bigr)
 =s^2(s-1)(s-2)x^{s-1}>0,
\]
so $h'''$ can vanish at most once. A function with two prescribed
double zeros and an additional zero would, by Rolle's theorem, force
at least two distinct zeros of $h'''$. Thus $h$ never vanishes outside
$a,b$. Its sign on $(0,a)$ is positive, and passing through a double
zero does not change that sign. Hence $h>0$ everywhere except at the
two interpolation nodes, as required.
\end{proof}
Apply this lemma with $s=1/5$, $a=1/32$, and $b=(37/40)^5$. Integrating $Q$ against $\nu_4$ uses only the three exact values in \eqref{eq:entropy}. The rational coefficient verification is recorded in Appendix~\ref{app:entropy}. It yields
\begin{equation}\label{eq:Cstar}
 \frac{A_{22/5}}{A_4}=\E_{\nu_4}X^{1/5}>\frac{101}{125}=:C_*.
\end{equation}

The resulting lower bound for $A_{22/5}$ will constrain both $S$ through \eqref{eq:peakbound} and $r$ through \eqref{eq:allspectral}.

The spectral estimates above control $r$, but the scalar expression
\eqref{eq:necessary} also depends on $S$ and $b_p$. As the exponent
changes, the effect of these factors must be controlled as well.
\subsubsection*{Scalar comparison across real exponents}
The following lemma carries a decay bound for $r$ into a decay
bound for $\cR_p$. It will be used here and in the final subrange,
where the degree-two spectral mass is retained explicitly.
\begin{lemma}[Scalar decay in the parameter]\label{lem:scalar-decay}
Let $2\leq S\leq7/2$ and $1/7\leq\beta\leq1/5$. For $0\leq t\leq1$ put $p=2+t$ and $x_0=\sqrt{\beta/(S-2/3)}$, and
\begin{equation}\label{eq:Cfunc}
 \mathcal C(S,\beta)=\frac14\log\frac{S-2/3}{S}-\frac1{18(S-2/3)}
 -\frac{\log S}{2(S-1)}+\frac56-\log2
 +\log\frac{F_3(x_0)}{F_2(x_0)}.
\end{equation}
Then $\mathcal C(S,\beta)-\beta/4$ increases with both $S$ and $\beta$; for $a\geq0$,
\begin{equation}\label{eq:scalar-est}
 \cR_{2+t}(\sqrt\beta e^{-at},S-1)
 \leq\cR_2(\sqrt\beta,S-1)\exp\{t(\mathcal C(S,\beta)-a)\}.
\end{equation}
In particular, writing $a_\beta=\beta/4-1/500$,
\begin{align}
 \cR_{2+t}(\sqrt\beta e^{-a_\beta t},S-1)
 &\leq e^{-t/10}\cR_2(\sqrt\beta,S-1)
 &&(2\leq S\leq7/2,\ 1/7\leq\beta\leq1/5),\label{eq:decay1}\\
 \cR_{2+t}(\sqrt\beta,S-1)
 &\leq e^{-2t/25}\cR_2(\sqrt\beta,S-1)
 &&(2\leq S\leq117/40,\ 1/7\leq\beta\leq3/16).
 \label{eq:decay2}
\end{align}
\end{lemma}
\begin{proof}
\par\medskip
\noindent\emph{Dependence on the exponent.}\quad
Write $r_a=\sqrt\beta e^{-at}$ and $x=r_a/b_p$. The interpolation
inequality for $F_p$ gives
\[
 \cR_p(r_a,S-1)\leq
 r_a K_p(S)\sqrt S\eta_p^{1-1/p}
 \frac{p(p+1)}{6\,2^t}F_2(x)^{1-t}F_3(x)^t.
\]

The quantity $b_p$ is nondecreasing in $p$:
\[
 \log b_p=\frac12\log S-\int_{1-S^{-1}}^1\frac{\dd u}{1+pu},
 \qquad
 \frac{\partial}{\partial p}\log b_p
 =\int_{1-S^{-1}}^1\frac{u}{(1+pu)^2}\dd u\geq0.
\]
Since $r_a$ is nonincreasing along $p=2+t$, it follows that
$x\leq x_0=\sqrt{\beta/(S-2/3)}$. Thus the $F_j$ terms can be
compared at a single fixed argument $x_0$.

Only the remaining scalar factors need to be differentiated. Set
\[
 V(p)=\log K_p(S)+(1-1/p)\log\eta_p
       +\log\frac{p(p+1)}{6\,2^{p-2}}.
\]
Its derivative decreases: indeed,
\[
 (\log K_p)''\leq p^{-2},\qquad
 \left(\log\frac{p(p+1)}{6\,2^{p-2}}\right)''
 =-p^{-2}-(p+1)^{-2},
\]
while
\[
 ((1-1/p)\log\eta_p)''
 =2\int_{1-1/S}^1\frac{u(1+u)}{(1+pu)^3}\dd u<\frac1{(p+1)^2}.
\]
Indeed, since $S\geq2$, the integration interval lies in
$[1/2,1]$. For $p\geq2$, the integrand
$u(1+u)/(1+pu)^3$ decreases on this interval; hence
\[
 2\int_{1-1/S}^1\frac{u(1+u)}{(1+pu)^3}\dd u
 \leq\frac{3/4}{(1+p/2)^3}<\frac1{(p+1)^2},
\]
where the last inequality is equivalent to $p^3+2>0$.
The three second-derivative bounds give $V''(p)<0$ and hence
$V'(p)\leq V'(2)$ for $p\geq2$.
We can therefore bound the change in $V$ by its derivative at $p=2$.
Since $F_2$ and $F_3$ are increasing on $[0,1)$, we may
replace their arguments by $x_0$ in the H\"older majorant above.
Using $p=2+t$, $r_a=\sqrt\beta e^{-at}$, and the definition of $V$,
we obtain the following comparison of function values:
\[
 \cR_{2+t}(r_a,S-1)
 \leq \sqrt\beta\,\sqrt S\,e^{-at+V(2+t)}
       F_2(x_0)^{1-t}F_3(x_0)^t.
\]
At $t=0$ the right-hand side equals
$\cR_2(\sqrt\beta,S-1)
 =\sqrt\beta\,\sqrt S\,e^{V(2)}F_2(x_0)$.
Consequently,
\[
 \frac{\cR_{2+t}(r_a,S-1)}{\cR_2(\sqrt\beta,S-1)}
 \leq \exp\!\left\{-at+V(2+t)-V(2)
          +t\log\frac{F_3(x_0)}{F_2(x_0)}\right\}.
\]
The concavity $V''<0$, established above, gives
$V(2+t)-V(2)\leq tV'(2)$.
A direct evaluation at $p=2$ shows that
\[
 V'(2)+\log\frac{F_3(x_0)}{F_2(x_0)}=\mathcal C(S,\beta).
\]
Substituting these two facts into the preceding inequality proves
\eqref{eq:scalar-est}. Thus the estimate follows from concavity of $V$, without comparing derivatives of $\log\cR_p$.

\par\medskip
\noindent\emph{Monotonicity in the auxiliary parameters.}\quad
For monotonicity in $\beta$, $S-2/3\leq17/6$ and $x_0<2/5$ imply
\[
 \partial_\beta\mathcal C
 =\frac{3(1+x_0)}{(S-2/3)(1+2x_0)F_3(x_0)}>\frac{250}{969}>\frac14.
\]
For the $S$ derivative, differentiate the ratio in \eqref{eq:Cfunc}
directly:
\[
 \partial_S\log\frac{F_3(x_0)}{F_2(x_0)}
 =-\frac{3x_0^2(1+x_0)}{(S-2/3)(1+2x_0)F_3(x_0)}.
\]
The quotient $x_0^2(1+x_0)/((1+2x_0)F_3(x_0))$ increases with
$x_0$ on the relevant interval. Hence the smallest possible
$\partial_S\mathcal C$ occurs at $\beta=1/5$. At $\beta=1/5$ the already defined $x_0$ satisfies
$S=(3+10x_0^2)/(15x_0^2)$. Using $1/4<x_0<2/5$ and
$\log S\geq2(S-1)/(S+1)$ gives
\[
 \begin{split}
 \frac1{x_0}\frac{\dd\mathcal C}{\dd x_0}
 &\leq-\frac5{3+10x_0^2}-\frac59
 -\frac{45}{(3+10x_0^2)(3+25x_0^2)}\\
 &\qquad+\frac{6(1+x_0)}{(1+2x_0)F_3(x_0)}
 <-\frac{335}{4347}<0.
 \end{split}
\]
Since $x_0$ decreases as $S$ increases, we obtain the asserted
monotonicity in $S$.

\par\medskip
\noindent\emph{Endpoint estimates for the parameter ranges.}\quad
For \eqref{eq:decay1}, evaluate the monotone quantity $\mathcal C-a_\beta$ at $S=7/2$, $\beta=1/5$. Here $x_0=\sqrt{6/85}<133/500$ and
$F_3(x_0)/F_2(x_0)<115/101$. The elementary logarithm bounds
\[
 \log(21/17)>211/1000,\quad \log(7/2)>313/250,\quad
 \log2>693/1000,\quad \log(115/101)<13/100
\]
then give $\mathcal C-a_\beta<-1/10$.
For \eqref{eq:decay2}, increase to $S=117/40$, $\beta=3/16$; now
$x_0=\sqrt{45/542}<15/52$,
$F_3(x_0)/F_2(x_0)<139/120$, and the corresponding bounds
\[
 \begin{gathered}
 \log(351/271)>129/500,\qquad \log(117/40)>134/125,\\
 \log2>693/1000,\qquad \log(139/120)<147/1000
 \end{gathered}
\]
show $\mathcal C<-2/25$. These inequalities follow from the positive arctanh expansion for logarithmic lower bounds and finite exponential-series lower bounds for the upper bounds. More precisely, for $y>1$ and $z=(y-1)/(y+1)$, the lower estimates use the truncated positive series
$\log y>2\sum_{j=0}^{n-1}z^{2j+1}/(2j+1)$, with $n=2,5,3,2,4$ for $y=21/17,7/2,2,351/271,117/40$, respectively. For the two upper estimates we use
$1+u+u^2/2+u^3/6>115/101$ at $u=13/100$ and
$1+u+u^2/2+u^3/6>139/120$ at $u=147/1000$. All comparisons are between rational numbers.
\end{proof}

\begin{proposition}[Exclusion on $22/5\leq q\leq5$]\label{prop:middle}
There are no noncharacter maximizers for $22/5\leq q\leq5$.
\end{proposition}
\begin{proof}
Write $q=4+t$ with $2/5\leq t\leq1$, and $p=2+t$.

\par\medskip
\noindent\emph{Control of the peak.}\quad
Convexity of $\log A_q$, using the secant from $4$ to $22/5$, and
\eqref{eq:Cstar} give
\begin{equation}\label{eq:middleconstant}
 \frac{A_{4+t}}{A_4}>C_*^{5t/2}.
\end{equation}
The peak bound \eqref{eq:peakbound} then forces
\[
 S<C_*^{-5}<\left(\frac{171}{100}\right)^2<\frac{117}{40}.
\]
The first comparison follows directly from $C_*^5(171/100)^2>1$.
\par\medskip
\noindent\emph{Control of the spectral radius.}\quad
By \eqref{eq:allspectral} and \eqref{eq:rbdef},
\begin{equation}\label{eq:middler}
 r^{2+t}\leq\frac{\overline G(t/2)}{5C_*^{5t/2}}.
\end{equation}
To prove $r<\sqrt3/4$, observe that the logarithm of the ratio
\[
 \frac{\overline G(t/2)}{5C_*^{5t/2}(\sqrt3/4)^{2+t}}
\]
is convex in $t$: $\overline G$ is a positive sum of exponentials,
while the other logarithmic terms are affine. A convex function lies
below the chord through its endpoints. Hence it suffices to show that
this ratio is less than one at $t=2/5$ and $t=1$. By \eqref{eq:Gend},
these are precisely the rational comparisons
\[
 (5C_*)^5(\sqrt3/4)^{12}>(27/50)^5,
 \qquad 22500C_*^5(\sqrt3/4)^6>49.
\]
Consequently
\begin{equation}\label{eq:middle-radius}
 r<R:=\frac{\sqrt3}{4}<\frac1{\sqrt5}.
\end{equation}
This explicit value of $R$ is crucial: $R^2=3/16$, which lies in the
parameter range allowed by the second estimate in \eqref{eq:decay2}.

It remains to exclude the two possibilities $1\leq S\leq2$ and
$2\leq S<117/40$. If $1\leq S\leq2$, Lemma~\ref{lem:smallpeaks} and \eqref{eq:middle-radius} immediately contradict \eqref{eq:necessary}. If $2\leq S<117/40$, use $r<R$, monotonicity in $r$, and \eqref{eq:decay2} with $\beta=R^2=3/16$:
\[
 \cR_{2+t}(r,S-1)\leq e^{-2t/25}
 \left(2R^2+R\sqrt{S-\tfrac23}\right).
\]
Since $S-2/3<271/120<(151/100)^2$,
\[
 2R^2+R\sqrt{S-\tfrac23}
 <2R^2+\frac{151R}{100}<\frac{103}{100},
\]
where $\sqrt3<26/15$ suffices. As $t\geq2/5$,
\[
 \cR_{2+t}(r,S-1)<\frac{103}{100}e^{-4/125}<1.
\]
This again contradicts \eqref{eq:necessary}, as required.
\end{proof}

\subsection{The range $4<q\leq22/5$: a coupled deficit and spectral contradiction}\label{subsec:lowrange}
The preceding estimates cover $22/5\leq q\leq5$. To treat $4<q\leq22/5$, separate bounds for the peak and the spectrum are not sufficient near $q=4$. The difficulty is concentrated in the second harmonic degree: it carries the largest weighted eigenvalue among the degrees allowed by \eqref{eq:peakdecomp}. If all admissible degrees are bounded by the same worst-case eigenvalue, the resulting estimate is too weak. We therefore retain the $\mathcal H_2$ mass as a common parameter in two estimates, one for the quartic deficit of $f$ and the other for the spectral radius $r$. The contradiction follows from requiring both estimates to hold for the same parameter. The final estimate must retain this common parameter rather than maximize the two bounds separately.

\subsubsection*{Keeping track of the degree-two mass}
Decompose the nonconstant part according to whether its harmonic degree is two or at least three. We shall measure their relative sizes by $\theta$, and use the same $\theta$ throughout the fourth-order and fractional-moment estimates.
For the unit $H\in\cH$ appearing in \eqref{eq:peakdecomp}, with
$P_0H=P_1H=0$, set
\begin{equation}\label{eq:bhd}
 b_H=\frac1{A_4}\int\phi^2(TH)^2\dd x,\quad
 h_H=\frac1{A_4}\int\phi(TH)^3\dd x,\quad
 d_H=\frac1{A_4}\int(TH)^4\dd x.
\end{equation}
Write
\begin{equation}\label{eq:theta}
 H=\sqrt\theta\,Y+\sqrt{1-\theta}\,G,
 \qquad\theta=\|P_2H\|_2^2\in[0,1],
\end{equation}
where $\|Y\|_2=\|G\|_2=1$, $Y\in\mathcal H_2$, and $P_0G=P_1G=P_2G=0$ (the endpoint cases are understood by continuity). Introduce
\begin{equation}\label{eq:thetas}
 \beta_\theta=\frac17+\frac{2\theta}{35},\qquad
 k_\theta=\frac47-\frac{6\theta}{35},\qquad
 H_\theta=\frac25-\frac{4\theta}{15}.
\end{equation}
The quartic spectrum \eqref{eq:spectra} implies $b_H\leq\beta_\theta$ and hence
$\kappa:=1-3b_H\geq k_\theta\in[2/5,4/7]$. Thus the single parameter
$\theta$ gives an explicit upper bound for the quadratic form $b_H$ and,
at the same time, a lower bound for the coefficient $\kappa$ in the
quartic deficit. This is why we do not replace $\theta$ by its worst
possible value at this stage.

We recall two mixed harmonic bounds, whose exact harmonic-block calculations are provided in Appendix~\ref{app:blocks}:
\begin{equation}\label{eq:mixedY}
 |h_Y|\leq\frac{2\sqrt5}{49},\quad
 \frac{\big|\int\phi(TY)^2TG\big|}{A_4}\leq\frac2{245},\quad
 \frac{\big|\int\phi(TY)(TG)^2\big|}{A_4}<\frac1{10}.
\end{equation}
For $G$ alone, the sharp quartic theorem applied to $G+\tau\one$ gives
a useful polynomial inequality. To see the coefficients, expand
\[
 \frac{\|T(G+\tau\one)\|_4^4}{A_4}
 =d_G+4\tau h_G+6\tau^2b_G+\tau^4,
\]
where the term cubic in $\tau$ vanishes by angular orthogonality, and compare
this with $(1+\tau^2)^2$, the fourth power of the $L^2$ norm.
After rearrangement, a quadratic polynomial in $\tau$ must be nonnegative
for every real $\tau$. Its discriminant is therefore nonpositive:
\begin{equation}\label{eq:discriminant}
 (2-6b_G)\tau^2-4h_G\tau+(1-d_G)\geq0,
 \qquad h_G^2\leq\frac{(1-3b_G)(1-d_G)}2.
\end{equation}
The other ingredient is Cauchy--Schwarz, in the form
\[
 \left(\int\phi(TG)^3\dd x\right)^2
 \leq\left(\int\phi^2(TG)^2\dd x\right)
      \left(\int(TG)^4\dd x\right),
 \qquad h_G^2\leq b_Gd_G.
\]
Combining this with \eqref{eq:discriminant} and eliminating $d_G$
leads to
\begin{equation}\label{eq:hG}
 h_G^2\leq\frac{b_G(1-3b_G)}{1-b_G}\leq\frac2{21}.
\end{equation}
Here $b_G\leq1/7$, so $1-b_G>0$; the second inequality is
exactly $(1-7b_G)(2-9b_G)\geq0$.

Expanding the cube of \eqref{eq:theta} produces four terms: the pure $Y$ term, the two mixed terms, and the pure $G$ term. Applying the corresponding estimate above to each term gives
\begin{equation}\label{eq:cubic}
 |h_H|\leq\frac{2\sqrt5}{49}\theta^{3/2}
 +\frac6{245}\theta\sqrt{1-\theta}
 +\frac3{10}\sqrt\theta(1-\theta)
 +\sqrt{\frac2{21}}(1-\theta)^{3/2}.
\end{equation}
\subsubsection*{The cubic term stays below the deficit coefficient}
The quartic deficit also contains the cubic quantity $h_H$. We need to
control it in terms of the same parameter $\theta$. The mixed estimates
\eqref{eq:mixedY} and the discriminant estimate for $G$ give the bound
\eqref{eq:cubic}; the following lemma shows that this bound stays a fixed
distance below the coefficient $H_\theta$.

\begin{lemma}[Uniform cubic-moment gap]\label{lem:cubicgap}
For every $H$ above,
\begin{equation}\label{eq:cubicgap}
 |h_H|<H_\theta-\frac1{480}.
\end{equation}
\end{lemma}
\begin{proof}
Set $x=\sqrt\theta$, $y=\sqrt{1-\theta}$. The strict rational bounds
$2\sqrt5/49<23/250$, $6/245<1/40$, and $\sqrt{2/21}<309/1000$ turn \eqref{eq:cubic} into
\[
 |h_H|<S(x,y):=\frac{23}{250}x^3+\frac1{40}x^2y
 +\frac3{10}xy^2+\frac{309}{1000}y^3.
\]
Since $(y-4/5)^2\geq0$, $y\leq41/40-5x^2/8$. Substitution, using $y^2=1-x^2$, yields
\[
 H_\theta-S(x,y)\geq P(x):=
 \frac{3331}{40000}-\frac3{10}x+\frac{26107}{120000}x^2
 +\frac{26}{125}x^3-\frac{71}{400}x^4.
\]
Writing $u=x-\frac12$, an exact expansion gives
\[
 P\Bigl(\frac12+u\Bigr)
 =\frac{617}{240000}-\frac{1823}{120000}u
 +\frac{31597}{120000}u^2-\frac{147}{1000}u^3-\frac{71}{400}u^4.
\]
For $|u|\leq\frac12$, the constant term is larger than $1/390$,
the absolute value of the linear coefficient is smaller than $1/65$,
and
\[
 \frac{31597}{120000}-\frac{147}{2000}
 -\frac{71}{1600}>\frac18.
\]
It follows that
\[
 P(x)>\frac1{390}-\frac{|u|}{65}+\frac{u^2}{8}
 \geq\frac1{390}-\frac2{65^2}>\frac1{480}.
\]
This proves the uniform gap.
\end{proof}
The uniform bound also applies at $\theta=0$ and $\theta=1$, so the quartic estimate covers both pure-mode cases.
\subsubsection*{A lower bound for the quartic deficit}
With this notation, the quartic deficit takes the form
\begin{equation}\label{eq:quarticdef}
 \delta_4(f):=1-\frac{\|Tf\|_4^4}{A_4}
 =2(1-3b_H)c^2\rho^2-4h_Hc\rho^3+(1-d_H)\rho^4.
\end{equation}
Applying the discriminant argument of \eqref{eq:discriminant}
with $H$ in place of $G$ gives $1-d_H\geq2h_H^2/\kappa$. Substituting into
\eqref{eq:quarticdef}, we can complete the square without dropping
the sign of the cubic term:
\[
 \begin{split}
 \delta_4(f)
 &\geq2\kappa c^2\rho^2-4h_Hc\rho^3
       +\frac{2h_H^2}{\kappa}\rho^4\\
 &=\frac{2\rho^2}{\kappa}(c\kappa-\rho h_H)^2.
 \end{split}
\]
Lemma~\ref{lem:cubicgap} controls $h_H$ by the fixed quantity
$H_\theta-1/480$, while $\kappa\geq k_\theta$. Because the squared
quantity may change sign, we retain its positive part before making
this replacement. It follows that
\begin{equation}\label{eq:deficitlower}
 \delta_4(f)\geq\frac{2\rho^2}{\kappa}
 \bigl[c\kappa-\rho h_H\bigr]^2
 \geq\frac{2\rho^2}{k_\theta}
 \left[c k_\theta-\rho\bigl(H_\theta-\tfrac1{480}\bigr)\right]_+^2.
\end{equation}
Here $[x]_+=\max\{x,0\}$. The monotonicity in $\kappa$ follows by writing the right side as
$2\rho^2\kappa[c-\rho(H_\theta-1/480)/\kappa]_+^2$.

\subsubsection*{The entropy bound for the same deficit}
We now bound the deficit from the other side. A maximizer at $q=4+t$ has $L^{4+t}$ norm at least that of a character, while its pointwise size is controlled by the peak $c$. With $0<t\leq2/5$, this comparison limits how large its quartic deficit can be. From \eqref{eq:entropy}, Jensen's inequality gives
\begin{equation}\label{eq:jensen}
 \frac{A_{4+t}}{A_4}=\E_{\nu_4}X^{t/2}\geq e^{-Mt},\qquad M=\frac{37}{60}.
\end{equation}
Since $A_{4+t}\leq\|Tf\|_{4+t}^{4+t}\leq c^tA_4(1-\delta_4(f))$, we deduce
\begin{equation}\label{eq:entropybudget}
 c\geq e^{-M},\qquad
 \delta_4(f)\leq1-e^{-t(M+\log c)}\leq t(M+\log c).
\end{equation}
The two inequalities retain the same $\theta$ as the spectral estimate below.

\subsubsection*{Spectral decay with the degree-two weight retained}
The spectral estimate provides the second restriction on the same $\theta$. We retain the separate rates of decay of the degree-two and higher-degree eigenvalues in the numerator of $r$.
We sharpen the two-node bounds \eqref{eq:PB} by checking
\[
 \frac{331}{10000}<x_-<\frac{207}{6250},\qquad
 \frac{57}{125}<x_+<\frac{457}{1000}.
\]
Direct substitution gives $P_B(207/6250)<0<P_B(457/1000)$; the signs themselves are all that is needed.
Together with $w>10/11$, elementary fifth-power comparisons imply
\[
 G(1/5)<\frac{10}{11}\frac{253}{500}+\frac1{11}\frac{107}{125}
 =\frac{1479}{2750}<e^{-31/50}.
\]
For the last strict inequality, the first three positive terms in the arctanh series establish
$\log(2750/1479)>31/50$. By log-convexity of $G$ with $G(0)=1$,
\begin{equation}\label{eq:expdegree2}
 \frac{\beta_{4+t,2}}{A_4}\leq\frac15e^{-31t/20}
 \qquad(0\leq t\leq2/5).
\end{equation}
For $k\geq3$, the higher-mode interpolation \eqref{eq:highmodes} gives
\begin{equation}\label{eq:exphigher}
 \frac{\beta_{4+t,k}}{A_4}\leq\frac17 30^{-t/2}
 \leq\frac17e^{-17t/10},
\end{equation}
using $\log30>17/5$, which follows from the arctanh series for $\log2$ and
$\log(16/15)<1/15-1/450+1/10125$.
By \eqref{eq:rbdef}, \eqref{eq:peakbound}, and the decomposition of $H$,
\begin{equation}\label{eq:coupledspectral}
 r^{2+t}\leq e^{Mt}
 \left(\frac{\theta}{5}e^{-31t/20}
 +\frac{1-\theta}{7}e^{-17t/10}\right).
\end{equation}
Because $\beta_\theta$ also appears in the deficit bound, it is essential that the next estimate retain its dependence on $\theta$.
\begin{lemma}[Radius retaining the quartic parameter]\label{lem:radius}
For $0\leq t\leq2/5$ and $\beta=\beta_\theta\in[1/7,1/5]$, the bound \eqref{eq:coupledspectral} implies
\begin{equation}\label{eq:radius}
 r\leq\sqrt\beta\,e^{-\beta t/4}.
\end{equation}
\end{lemma}
\begin{proof}
Recall that $\beta=1/7+2\theta/35$. The right side of
\eqref{eq:coupledspectral} can therefore be written entirely in
terms of $\beta$:
\[
 e^{Mt}\beta e^{-17t/10}
 \left[1+\frac{7\beta-1}{2\beta}
                 (e^{3t/20}-1)\right].
\]
For $0\leq t\leq2/5$, the exponential series gives
$e^{3t/20}-1\leq(31/200)t$. Using $\log(1+y)\leq y$,
\[
 \log r\leq
 \frac{\log\beta+
  \bigl(M-\frac{463}{400}-\frac{31}{400\beta}\bigr)t}{2+t}.
\]
In view of $M=37/60$ and $t\leq2/5$, the desired inequality
\eqref{eq:radius} follows if
\[
 \frac{649}{1200}+\frac{31}{400\beta}
 +\frac12\log\beta-\frac35\beta>0.
\]
The left side is decreasing in $\beta$: multiplying its derivative by
$\beta^2$ gives
$-31/400+\beta/2-3\beta^2/5<0$ on $[1/7,1/5]$.
At $\beta=1/5$ it equals
$97/120-\tfrac12\log5>1/1200$, since $\log5<323/200$.
This proves \eqref{eq:radius}.
\end{proof}
The last elimination is simplest in terms of $k=k_\theta$ from
\eqref{eq:thetas} rather than $\theta$: the three quantities satisfy
$\beta_\theta=(1-k)/3$ and
$H_\theta=14k/9-22/45$. Thus this change is not an additional
independent parameter; it expresses the two coupled estimates in the
same variable.

\subsubsection*{The coupled restrictions on the peak and the spectrum}
The quartic estimate and the spectral estimate now constrain the same two parameters, but in different forms. It is useful to express both in terms of the peak ratio $v=\rho/c$ and the parameter $k=k_\theta$. Put $S=1+v^2$ and abbreviate $k=k_\theta\in[2/5,4/7]$. Define
\begin{equation}\label{eq:CvDk}
 e(v)=\frac{37}{60}-\frac12\log(1+v^2),\quad
 C(v)=\frac{(1+v^2)^2}{2v^2},\quad
 D(k)=\frac{14}{9}k-\frac{707}{1440}.
\end{equation}
Recall that $c=(1+v^2)^{-1/2}$, $\rho=vc$ and
$M+\log c=e(v)$, with $M=37/60$ from \eqref{eq:jensen}.
Substituting these identities into \eqref{eq:deficitlower} and comparing
with the upper bound \eqref{eq:entropybudget}, then cancelling the common
positive factor $2v^2c^4/k$, gives the first necessary condition
\begin{equation}\label{eq:budget2}
 [k-vD(k)]_+^2\leq t\,C(v)e(v)k,\qquad e(v)\geq0.
\end{equation}
The right-hand side of \eqref{eq:budget2} tends to zero with $t$. Hence, near the endpoint, the parameters must approach the relation $k=vD(k)$ unless the positive part vanishes. That correlation is the information which a separate estimate on $v$ or $k$ would erase. Independently, the spectral radius bound controls the quantity in the necessary condition. Indeed, $\beta_\theta=(1-k)/3$ and
\[
 r\leq\sqrt{\beta_\theta}\,e^{-\beta_\theta t/4}
 \leq\sqrt{\beta_\theta}\,e^{-a_{\beta_\theta}t},
 \qquad a_{\beta_\theta}=\frac{\beta_\theta}{4}-\frac1{500}.
\]
Thus, by monotonicity of $\cR_p$ in $r$ and \eqref{eq:decay1}, whenever
$2\leq S\leq7/2$ the second necessary bound is
\begin{equation}\label{eq:R2decay}
 \cR_{2+t}(r,v^2)\leq e^{-t/10}\,\mathcal B(v,k),\qquad
 \mathcal B(v,k):=\cR_2\!\left(\sqrt{\frac{1-k}{3}},v^2\right)
 =\frac{2(1-k)+\sqrt{(1-k)(3v^2+1)}}3.
\end{equation}
The two necessary restrictions now involve the same $k$: rule \eqref{eq:budget2} limits the pairs $(v,k)$ allowed by the quartic deficit, while \eqref{eq:R2decay} bounds the scalar quantity to be contradicted. We use the first restriction to control the second, rather than maximizing separately over $v$ and $k$.

For $v$ close to its largest admissible value, the entropy
factor $e(v)$ is small. The next estimate makes that observation
quantitative and is used only in the final two-parameter elimination.
\begin{lemma}[Entropy budget near its limiting peak]\label{lem:ebudget}
Let $U=39/25$ and $\gamma_*=1073/975$. If $v\geq6/5$ and $e(v)\geq0$, then
\begin{equation}\label{eq:ebud}
 v<U,\qquad C(v)e(v)<\gamma_*(U-v).
\end{equation}
\end{lemma}
\begin{proof}
The positive arctanh series gives the explicit lower bound
\[
 \log(1+U^2)>
 2\sum_{j=0}^{5}\frac{(1521/2771)^{2j+1}}{2j+1}>\frac{37}{30}.
\]
Thus $e(v)\geq0$ implies
$1+v^2\leq e^{37/30}<1+U^2$, so $v<U$. By
$\log(b/a)\leq(b-a)/\sqrt{ab}$ for $b\geq a>0$,
\[
 C(v)e(v)
 <\frac{(1+v^2)^{3/2}(U+v)}{4v^2\sqrt{1+U^2}}(U-v).
\]
The logarithmic derivative of the factor multiplying $U-v$ has the sign of
$U(v^2-2)+v(2v^2-1)$. This polynomial is positive at $v=6/5$ and strictly increasing for $v\geq6/5$, since its derivative is $2Uv+6v^2-1>0$. Therefore the factor is increasing on $[6/5,U]$ and is bounded by its endpoint value $(1+U^2)/(2U)=\gamma_*$.
\end{proof}
The remaining obstruction is that the quartic budget must hold at
the same $(v,k)$ as the endpoint majorant $\mathcal B$. Writing
$u=U-v$ makes this coupling a quadratic inequality in $u$.
The matrix argument below eliminates the pair without maximizing the
two estimates independently.
\begin{proposition}[Positive-quadratic-form elimination]\label{prop:quadratic}
Suppose $0\leq t\leq2/5$, $2/5\leq k\leq4/7$, $v\geq6/5$, $e(v)\geq0$, and \eqref{eq:budget2} holds. Then
\begin{equation}\label{eq:R2upper}
 \mathcal B(v,k)<\frac{73}{75}+\frac t6\leq1+\frac t{10}.
\end{equation}
\end{proposition}
\begin{proof}
By Lemma~\ref{lem:ebudget} put $u=U-v>0$. Inequality \eqref{eq:budget2} implies
\begin{equation}\label{eq:Qbudget}
 [Q(k)+D(k)u]_+^2\leq t\gamma_*ku,\qquad
 Q(k):=k-UD(k)=\frac{9191}{12000}-\frac{107}{75}k.
\end{equation}
On $v\in[6/5,U]$, the affine function $(128v+31)/80$ majorizes $\sqrt{3v^2+1}$. Indeed, the difference is concave in $v$, and direct squaring verifies its positivity at both endpoints. The tangent-line bound at $k=7/16$ is
$\sqrt{1-k}\leq(25-16k)/24$. Substituting these two elementary
majorants directly into the explicit expression for $\cR_2$ yields
\begin{equation}\label{eq:affinebound}
 \mathcal B(v,k)\leq
 \frac{640v+923}{1152}-\frac{128v+271}{360}\,k.
\end{equation}
The right side is affine in $v$, so
\begin{equation}\label{eq:Delta}
 \mathcal B(v,k)\leq\frac{73}{75}+\Delta(k)-\ell(k)u,
 \quad \Delta(k)=\frac{20003}{28800}-\frac{11767}{9000}k,
 \quad\ell(k)=\frac{25-16k}{45}>0.
\end{equation}
If $\Delta\leq0$, the result follows at once. Assume $\Delta>0$, and put
$m(k)=9/8+(k-2/5)/3$. For $h=k-2/5$, expanding the left-hand side and bounding its coefficients gives
\begin{equation}\label{eq:QDelta}
 Q-m\Delta>\frac1{500}-\frac h{60}+\frac{2h^2}{5}
 \geq\frac1{500}-\frac1{5760}>0.
\end{equation}
In particular, $Q>m\Delta>0$ and $D(k)\geq21/160>0$.
A second direct expansion, using $0\leq h\leq6/35$, gives the uniform comparison
\begin{equation}\label{eq:positive-poly}
 \frac23m(D+m\ell)-\gamma_*k
 >\frac1{150}+\frac h{1000}+\frac{h^2}{6}>0.
\end{equation}
In particular, the positive margin does not degenerate at either endpoint of the $k$-interval. Inequalities \eqref{eq:QDelta} and \eqref{eq:positive-poly} are arranged to verify the two principal-minor conditions for the matrix below.
Consider the symmetric matrix
\begin{equation}\label{eq:matrix}
 \mathsf A=
 \begin{pmatrix}
 Q^2/6&QD/6-\gamma_*k\Delta/2\\
 QD/6-\gamma_*k\Delta/2&D^2/6+\gamma_*k\ell
 \end{pmatrix}.
\end{equation}
Its upper-left entry $Q^2/6$ is positive because $Q>m\Delta>0$.
To verify the determinant, multiply out the two diagonal entries and
subtract the square of the off-diagonal entry:
\[
 \begin{split}
 \det\mathsf A
 &=\frac{Q^2D^2}{36}+\frac{\gamma_*kQ^2\ell}{6}
   -\left(\frac{Q^2D^2}{36}
     -\frac{\gamma_*kQD\Delta}{6}
     +\frac{\gamma_*^2k^2\Delta^2}{4}\right)\\
 &=\frac{\gamma_*k}{4}
 \left(\frac23Q(D\Delta+Q\ell)-\gamma_*k\Delta^2\right).
 \end{split}
\]
Because $D,\Delta,\ell>0$, the inequality $Q>m\Delta$ gives
\[
 Q(D\Delta+Q\ell)>m\Delta^2(D+m\ell).
\]
The second polynomial comparison \eqref{eq:positive-poly} says that
$(2/3)m(D+m\ell)>\gamma_*k$; hence the bracket in the determinant
is strictly positive. Sylvester's criterion now makes $\mathsf A$
positive definite.
Positive definiteness of this matrix compares the quartic
budget directly with the possible excess $\Delta-\ell u$, uniformly
in $k$. Evaluating its quadratic form at $(1,u)$ gives
\[
 \frac16(Q+Du)^2-\gamma_*ku(\Delta-\ell u)>0.
\]
The quantity $Q+Du$ is positive, and \eqref{eq:Qbudget} yields
$\Delta-\ell u<(Q+Du)^2/(6\gamma_*ku)\leq t/6$.
Now \eqref{eq:Delta} proves
$\mathcal B(v,k)<73/75+t/6$. Finally
$73/75+t/6-(1+t/10)=-2/75+t/15\leq0$.
\end{proof}
The matrix estimate handles $v\geq6/5$. The remaining values of $v$ are covered by the scalar comparison lemmas, so we can now conclude the argument for $4<q\leq22/5$.
\begin{proposition}[Exclusion on $4<q\leq22/5$]\label{prop:low}
No noncharacter maximizer exists for any real $4<q\leq22/5$.
\end{proposition}
\begin{proof}
Let $q=4+t$, $0<t\leq2/5$, and suppose $f$ is such a maximizer. Its parameters satisfy \eqref{eq:budget2}.

\medskip
\noindent\textbf{Small peaks: $S\leq2$.}\par
If $S\leq2$, Lemma~\ref{lem:radius} gives $r\leq1/\sqrt5$ and
Lemma~\ref{lem:smallpeaks} contradicts \eqref{eq:necessary}.

\medskip
\noindent\textbf{Moderate peaks: $2<S\leq61/25$.}\par
Now assume $S>2$, hence $v>1$. If $1<v\leq6/5$, maximizing the explicit expression in \eqref{eq:R2decay} over $k\in[2/5,4/7]$ and $v\in[1,6/5]$ gives
\[
 \mathcal B(v,k)\leq\frac25+\sqrt{\frac{133}{375}}<1.
\]
Then $S\leq61/25<7/2$, and \eqref{eq:R2decay} proves $\cR_{2+t}<1$.

\medskip
\noindent\textbf{Large peaks: $v\geq6/5$.}\par
If $v\geq6/5$, Lemma~\ref{lem:ebudget} implies
$S<2146/625<7/2$, and Proposition~\ref{prop:quadratic} together with \eqref{eq:R2decay} gives
\[
 \cR_{2+t}(r,v^2)
 <e^{-t/10}\left(\frac{73}{75}+\frac t6\right)
 \leq e^{-t/10}(1+t/10)\leq1.
\]
In every case we contradict \eqref{eq:necessary}. Existence follows from Theorem~\ref{thm:exist}, and the reduction \eqref{eq:realreduction} gives the complex extremizers.
\end{proof}

The three subranges cover $4<q\leq6$; Proposition~\ref{prop:base} then extends the conclusion to higher exponents.

\begin{proof}[Completion of Theorem~\ref{thm:main}]
Theorem~\ref{thm:foschi} supplies the already known quartic endpoint.
The arguments of Subsections~\ref{subsec:highrange}, \ref{subsec:middle}
and \ref{subsec:lowrange} establish the sharp inequality and complete
classification on $(4,6]$. Proposition~\ref{prop:base} propagates both
assertions to every finite real $q>4$.

At $q=\infty$, since $\mu$ has mass one, Cauchy--Schwarz gives
\[
 |Tf(x)|=\left|\int_{\Sph}f(\omega)e^{-ix\cdot\omega}
                 \dd\mu(\omega)\right|
 \leq\|f\|_{L^2(\mu)}\qquad(x\in\R^3).
\]
The constant function attains equality at $x=0$, so
$\|T\|_{2\to\infty}=1$. Conversely, suppose that a nonzero $f$
attains this operator norm. As established in Section~\ref{sec:prelim},
$Tf\in C_0(\R^3)$ and $Tf$ is not identically zero. Hence there is
$a\in\R^3$ such that
\[
 |Tf(a)|=\|Tf\|_\infty=\|f\|_{L^2(\mu)}.
\]
Equality in Cauchy--Schwarz implies that
$f(\omega)=ce^{ia\cdot\omega}$ almost everywhere on $\Sph$, for
some $c\in\C\setminus\{0\}$. Conversely, every such function attains
the norm, since $Tf(a)=c$ and $\|f\|_2=|c|$.

Finally, $\widehat{f\sigma}=4\pi Tf$ and
$\|f\|_{L^2(\sigma)}=\sqrt{4\pi}\|f\|_{L^2(\mu)}$ give the
area-measure constants in \eqref{eq:area} and at $q=\infty$.
\end{proof}

\appendix
\section{Exact convolution formulas and entropy identities}\label{app:convol}
The radial convolution identities below yield the moments $A_{2m}$
defined in \eqref{eq:phi} and the harmonic eigenvalues $\beta_{q,k}$
introduced in Subsection~\ref{subsec:common-spectral}. They also fix the
normalizations in the logarithmic-moment calculation.

Under the normalized spherical measure $\mu=\sigma/(4\pi)$,
the first coordinate has density $\rho(t)=\frac12\one_{[-1,1]}(t)$. Its convolution powers are piecewise polynomials; for the moment and spectral calculations below, we only need their expressions on $0\leq r\leq2$. For even $n$, the formula is
\begin{equation}\label{eq:gn}
 g_n(r)=\frac1{2^n(n-1)!}\sum_{j=0}^{n/2}(-1)^j\binom nj(r+n-2j)^{n-1}.
\end{equation}
Indeed, write
$\rho(t)=\frac12\bigl(\one_{t>-1}-\one_{t>1}\bigr)$.
In the $n$-fold convolution expand the $n$ factors by the binomial formula.
A term selecting $j$ copies of $\one_{t>1}$ and $n-j$ copies of
$\one_{t>-1}$ is the $(n-1)$st positive-part monomial, translated
by $n-2j$. This follows by induction from
\[
 \frac{t_+^a}{a!}*\frac{t_+^b}{b!}
 =\frac{t_+^{a+b+1}}{(a+b+1)!},
\]
where the identity is the beta integral on $0\leq s\leq t$.
The coefficients are $(-1)^j\binom nj/2^n$. For $0\leq r\leq2$,
all the terms with $j>n/2$ vanish, giving \eqref{eq:gn}.
The restriction to a fixed radial interval is important here: the full
convolution is piecewise polynomial, while these displayed terms are
sufficient for the values and derivatives at the origin used below.

Let $M_n(r)$ be the radial density of $\mu^{*n}$. Integration in the transverse coordinates gives
\[
 g_n(t)=2\pi\int_{|t|}^\infty M_n(r)r\dd r,
 \qquad M_n(r)=-\frac{g_n'(r)}{2\pi r}.
\]
Since $\phi^n\in L^1(\R^3)$ for $n\geq4$, Fourier inversion gives
$\int\phi^n\dd x=(2\pi)^3M_n(0)$ with $M_n(0)=-g_n''(0)/(2\pi)$. In particular,
\begin{equation}\label{eq:AnTable}
\begin{array}{c|cccc}
 n&4&6&8&10\\\hline
 -g_n''(0)&1/4&1/8&1/12&35/576\\
 A_n/\pi^2&1&1/2&1/3&35/144
\end{array}
\end{equation}
The polynomials needed for these derivatives are
\begin{align*}
 g_2(r)&=\tfrac12-\tfrac14r,\\
 g_4(r)&=\tfrac13-\tfrac18r^2+\tfrac1{32}r^3,\\
 g_6(r)&=\tfrac{11}{40}-\tfrac1{16}r^2+\tfrac1{128}r^4-\tfrac1{768}r^5,\\
 g_8(r)&=\tfrac{151}{630}-\tfrac1{24}r^2+\tfrac1{288}r^4-\tfrac1{4608}r^6+\tfrac1{36864}r^7,\\
 g_{10}(r)&=\tfrac{15619}{72576}-\tfrac{35}{1152}r^2+\tfrac{19}{9216}r^4
 -\tfrac5{55296}r^6+\tfrac1{294912}r^8-\tfrac1{2949120}r^9.
\end{align*}

For $n\geq4$, the kernel of $T^*\phi^nT$ is $(2\pi)^3M_n(|\omega-\eta|)$. In the case $n=2$, regularize with a Gaussian in Fourier space; $M_2(y)=(8\pi|y|)^{-1}\one_{|y|<2}$, and the regularized Newtonian potential satisfies
\[
 G_\varepsilon*|\cdot|^{-1}(y)=|y|^{-1}\operatorname{erf}(|y|/(2\sqrt\varepsilon))\leq |y|^{-1}.
\]
Dominated convergence on the sphere, followed by density and the quartic bound, justifies this singular-kernel limit. For $k\geq1$ the corresponding zonal kernels are $\pi^2/r$ and $\pi^2(1-3r/8)$ for $\beta_{4,k}$ and $\beta_{6,k}$, respectively. Funk--Hecke applied to the general chord power gives
\begin{equation}\label{eq:chord}
 \frac12\int_{-1}^1(2-2t)^{a/2}\mathsf P_k(t)\dd t
 =\frac{2^{a+1}(-a/2)_k}{(a+2)(a/2+2)_k},\qquad a>-2,
\end{equation}
where $(z)_k$ is the rising factorial. For completeness, substitute $r=(2-2t)^{1/2}$ in the Funk--Hecke
integral and apply Rodrigues' formula
\[
 \mathsf P_k(t)=\frac1{2^kk!}\frac{\dd^k}{\dd t^k}(t^2-1)^k.
\]
After $k$ integrations by parts the differentiation falls on
$(1-t)^{a/2}$, and the remaining integral is a beta integral.
The resulting products simplify to
$2^{a+1}(-a/2)_k/((a+2)(a/2+2)_k)$.
The initial calculation is valid where the integrations by parts have
vanishing endpoint terms; both sides are analytic for $a>-2$
after continuation, proving \eqref{eq:chord} in the stated range.
At $a=-1$, the formula gives the Newtonian eigenvalue
$(2k+1)^{-1}$. At $a=1$ it gives the degree-$k$ eigenvalues
of the chord-length kernel. Substituting these values into the
explicit kernels for $T^*\phi^2T$ and $T^*\phi^4T$ yields
\eqref{eq:spectra} for every required $k$.

For a unit degree-two harmonic, $\mathsf P_2(1-r^2/2)=1-\frac32r^2+\frac38r^4$. Using $M_n=-g_n'/(2\pi r)$,
\begin{equation}\label{eq:b2gn}
 \frac{\beta_{n+2,2}}{\pi^2}=-2\int_0^2g_n'(r)
 \left(1-\frac32r^2+\frac38r^4\right)\dd r.
\end{equation}
For $n=2,4,6,8$ the values are respectively $1/5$, $1/70$, $1/252$ and $41/23760$. Dividing the last three by $1/5$ proves \eqref{eq:three-moments}.

\subsection*{Logarithmic differentiation of the fourth moment}\label{app:entropy}
Recall the fourth-power probability measure
$\dd\nu_4=A_4^{-1}\phi^4\dd x$ and $X=\phi^2$ introduced
before \eqref{eq:entropy}. We compute the average of $-\log X$
under $\nu_4$ by separating the contributions of $\log|\sin r|$
and $\log r$, both obtained by differentiating parameter-dependent integrals.
For $p>1$, periodization and the partial-fraction identity for $\csc^2$ give
\begin{equation}\label{eq:period}
 \int_0^\infty\frac{|\sin r|^p}{r^2}\dd r
 =\frac12\int_0^\pi\sin^{p-2}u\dd u
 =\frac{\sqrt\pi}{2}\frac{\Gamma((p-1)/2)}{\Gamma(p/2)}.
\end{equation}
To identify the logarithmic moment of the sine, differentiate the
logarithm of the right side of \eqref{eq:period}. The formula
$\psi(z)=\Gamma'(z)/\Gamma(z)$ gives
\[
 \left.\frac{\dd}{\dd p}\log
 \frac{\Gamma((p-1)/2)}{\Gamma(p/2)}\right|_{p=4}
 =\frac12\bigl(\psi(3/2)-\psi(2)\bigr)
 =\frac12-\log2.
\]
Because the normalized measure $\nu_4$ has radial density proportional
to $\sin^4r/r^2$, this is precisely
$\E_{\nu_4}\log|\sin r|=1/2-\log2$.
To compute the logarithmic radial moment, one uses, for $0<\nu<2$,
\[
 \int_0^\infty(1-\cos(ar))r^{-1-\nu}\dd r
 =\frac{\pi a^\nu}{2\Gamma(1+\nu)\sin(\pi\nu/2)}.
\]
Indeed, initially for $0<\nu<1$ this follows from an exponential cutoff and the Gamma integral; cancellation at the origin extends it to $0<\nu<2$. Since
$\sin^4r=\frac12(1-\cos2r)-\frac18(1-\cos4r)$, we get
\begin{equation}\label{eq:mellin}
 \int_0^\infty\sin^4r\,r^{-1-\nu}\dd r
 =\frac{\pi(4\cdot2^\nu-4^\nu)}{16\Gamma(1+\nu)\sin(\pi\nu/2)}.
\end{equation}
Differentiating \eqref{eq:mellin} at $\nu=1$ yields the logarithmic
radial moment. Indeed, if $M(\nu)$ denotes the integral on its left,
then $-M'(1)/M(1)$ is the $\nu_4$-average of $\log r$.
The differentiation of the numerator $4\cdot2^\nu-4^\nu$,
together with $\psi(2)=1-\gamma$ and the vanishing of the derivative
of $\log\sin(\pi\nu/2)$ at $\nu=1$, gives
\[
 \E_{\nu_4}\log r=1-\gamma.
\]
Hence
\[
 \E_{\nu_4}(-\log X)
 =-2\E_{\nu_4}\log|\sin r|+2\E_{\nu_4}\log r
 =1+2\log2-2\gamma.
\]
Differentiation under the integrals is justified by integrable logarithmic majorants, both at the sine zeros and at $r=0,\infty$. For the explicit upper bound in \eqref{eq:entropy}, write $H_{16}=\sum_{j=1}^{16}j^{-1}$. The trapezoidal integral comparison gives
$\gamma>H_{16}-4\log2-1/32$. Four terms of the positive arctanh series together with its geometric remainder give $\log2<1733/2500$. These yield
\[
 \E_{\nu_4}(-\log X)
 <\frac{17}{16}+10\frac{1733}{2500}-2H_{16}
 <\frac{37}{30}.
\]

\subsection*{Exact arithmetic for the logarithmic minorant}
The proof of Lemma~\ref{lem:logHermite} reduces the fractional-moment estimate to integrating its Hermite minorant against $\nu_4$. We record the interpolation coefficients and the rational bounds used in that calculation, so that the numerical constant in \eqref{eq:Cstar} can be checked directly.
For Lemma~\ref{lem:logHermite}, take $s=1/5$, $a=1/32$, $b=(37/40)^5$ and set
\begin{align*}
 u&=\frac1{10a},& v&=\frac{37}{200b},&
 C_0&=\frac{v-u}{2(b-a)},& B_0&=u-2aC_0,\\
 J&=\frac{b-a}{2}(u+v)-\frac{17}{40},&
 D_0&=\frac12\left(\frac ba-\frac ab\right),&
 D&=\frac{J}{D_0-5\log(37/20)}.
\end{align*}
The interpolation equations give
$C=C_0+D/(2ab)$, $B=B_0-D(1/a+1/b)$, and
$A=a^{1/5}-Ba-Ca^2-D\log a$. Thus
\[
 \E_{\nu_4}Q(X)=M_0+D\left(K_0+5\log2-\E_{\nu_4}(-\log X)\right),
\]
where
\[
 M_0=\frac12+B_0(1/2-a)+C_0(1/3-a^2)>\frac{6567}{5000},
\]
\[
 K_0=-\left(\frac1a+\frac1b\right)(1/2-a)+\frac{1/3-a^2}{2ab}.
\]
Using the rational logarithmic estimates $\log2>693/1000$,
$\log(37/20)<77/125$, one checks
\[
 0<D<\frac{2253}{25000},\qquad
 K_0+5\log2-\frac{37}{30}>-\frac{701}{125}.
\]
The latter logarithm bound follows by taking $z=17/57$ in the positive arctanh series:
\[
 \log(37/20)<2\sum_{j=0}^2\frac{z^{2j+1}}{2j+1}
 +\frac{2z^7}{7(1-z^2)}<\frac{77}{125}.
\]
For the bound on $D$, substitute $a=1/32$ and $b=(37/40)^5$ into the displayed formulas for $J$ and $D_0$, then use $\log(37/20)<77/125$. Clearing the positive denominator yields $D<2253/25000$ directly; as required for the final estimate.
Since the sign of the bracket is negative and $D>0$, these estimates yield
\[
 \E_{\nu_4}X^{1/5}>\E_{\nu_4}Q(X)
 >\frac{6567}{5000}-\frac{701}{125}\frac{2253}{25000}
 >\frac{101}{125},
\]
which proves \eqref{eq:Cstar}.

\section{Quadratic-harmonic couplings and the mixed operator}\label{app:blocks}
We establish the mixed harmonic estimates used near $q=4$ in
\eqref{eq:mixedY}. These calculations are independent of the target
exponent and exploit the banded action of multiplication by a
degree-two spherical harmonic. We first state the three bounds, then compute the harmonic blocks and estimate their action on the higher-degree subspace.
\begin{proposition}\label{prop:blocks}
Let $Y$ be a real degree-two spherical harmonic with $\|Y\|_{L^2(\mu)}=1$. Then
\begin{equation}\label{eq:block-bound1}
 \frac{\int\phi^2(TY)^2\dd x}{A_4}=\frac15,
 \qquad \frac{\big|\int\phi(TY)^3\dd x\big|}{A_4}\leq\frac{2\sqrt5}{49}.
\end{equation}
For $G\in\cH$ (with $\cH$ as in \eqref{eq:H}) satisfying
$P_0G=P_1G=P_2G=0$,
\begin{equation}\label{eq:block-bound2}
 \frac{\big|\int\phi(TY)^2TG\dd x\big|}{A_4}\leq\frac2{245}\|G\|_2.
\end{equation}
Moreover, for any complex such $G$,
\begin{equation}\label{eq:block-bound3}
 \frac{\big|\int\phi(TY)|TG|^2\dd x\big|}{A_4}
 \leq M_*\|G\|_2^2,
 \qquad M_*:=\frac{617\sqrt5}{13860}<\frac1{10}.
\end{equation}
\end{proposition}

\subsection*{Spherical pairings and the degree-two operator}
A real degree-two spherical harmonic has the form $Y(\omega)=\omega^\top A\omega$ with $A$ symmetric and trace-free. This representation lets us compute its moments and its interaction with other harmonic degrees using matrix traces. The basic integration formula is the spherical version of the Gaussian pairing identity.
For vectors $v_1,\ldots,v_{2m}\in\R^d$, the Gaussian pairing identity of Isserlis \cite{Isserlis}, in its spherical form \cite{Vignat}, reads
\begin{equation}\label{eq:pairing}
 \int_{\mathbb S^{d-1}}\prod_{j=1}^{2m}(v_j\cdot\omega)\dd\mu_d(\omega)
 =\frac{\displaystyle\sum_{\pi\in\mathcal P_2(2m)}\prod_{\{i,j\}\in\pi}(v_i\cdot v_j)}
 {d(d+2)\cdots(d+2m-2)}.
\end{equation}
It follows by writing a standard Gaussian as $R\omega$ and dividing Wick's formula by
$\E R^{2m}=d(d+2)\cdots(d+2m-2)$.

Write $Y(\omega)=\omega^\top A\omega$, where $A$ is a real
symmetric trace-free $3\times3$ matrix. Its normalization and
\eqref{eq:pairing} give
\[
 \int Y^2\dd\mu=\frac2{15}\tr(A^2)=1,
 \qquad \tr(A^2)=\frac{15}{2}.
\]
For a trace-free $3\times3$ matrix, $\tr(A^4)=\tfrac12(\tr A^2)^2$.
Applying the same pairing formula once more gives
\[
 \int Y^4\dd\mu
 =\frac{12(\tr A^2)^2+48\tr A^4}{945}=\frac{15}{7}.
\]
For a trace-free symmetric matrix $C$ it also gives
$\int Y^2(\omega^\top C\omega)\dd\mu=8\tr(A^2C)/105$.
Thus
\begin{equation}\label{eq:P2Y}
 P_2(Y^2)=\frac47\,\omega^\top\Bigl(A^2-\frac52 I\Bigr)\omega,
 \qquad\|P_2(Y^2)\|_2^2=\frac{20}{49}.
\end{equation}
Since $P_0(Y^2)=1$ and $Y^2$ has only degrees $0,2,4$,
\begin{equation}\label{eq:Y4P4}
 \|P_4(Y^2)\|_2^2=\frac{15}{7}-1-\frac{20}{49}=\frac{36}{49},
 \qquad\left|\int Y^3\dd\mu\right|
 =|\ip Y{P_2(Y^2)}|\leq\frac{2\sqrt5}{7}.
\end{equation}

Introduce the bounded self-adjoint operator
\begin{equation}\label{eq:KY}
 K_Y:=A_4^{-1}T^*(\phi\,TY)T.
\end{equation}
Conditioning on the sum of two unit vectors and using
$\mu*\mu(r\theta)=(8\pi r)^{-1}\one_{0<r<2}$ gives its kernel, with $r=|\omega-\eta|$,
\begin{equation}\label{eq:KYkernel}
 K_Y(\omega,\eta)
 =\left(\frac{3}{8r}-\frac1{2r^3}\right)
 (\omega-\eta)^\top A(\omega-\eta).
\end{equation}
The singularity is integrable; Appendix~\ref{app:convol} gives the Gaussian regularization needed to justify the kernel identity on general $L^2$ inputs.

In the commutator representation below, $D$ is the spectral
operator $D|_{\mathcal H_k}=2k+1$, and $M_i$ denotes multiplication
by $\omega_i$. The $r^{-1}$ kernel has eigenvalue $(2k+1)^{-1}$; the eigenvalue differences for the renormalized $r^{-3}$ kernel are $-k/2$:
\[
 \frac12\int_{-1}^1\frac{\mathsf P_k(t)-1}{(2-2t)^{3/2}}\dd t=-\frac{k}{2}.
\]
The latter follows by integration by parts and the Legendre derivative expansion
\[
 \mathsf P_k'=\sum_{j<k,\ k-j\text{ odd}}(2j+1)\mathsf P_j.
\]
An additive scalar disappears under double commutation. Formula \eqref{eq:KYkernel} therefore becomes
\begin{equation}\label{eq:commutator}
 K_Y=\frac38\sum_{i,j=1}^3 A_{ij}[M_i,[M_j,D^{-1}]]
 +\frac18\sum_{i,j=1}^3 A_{ij}[M_i,[M_j,D]].
\end{equation}
Multiplication by a coordinate sends $\mathcal H_k$ into $\mathcal H_{k-1}\oplus\mathcal H_{k+1}$. In the contracted products the coefficients of the lower and upper component are respectively $(n+2)/(2n)$ and $(n-2)/(2n)$, where $n=2k+1$. Inserting these into \eqref{eq:commutator}, or using the spherical product rule directly, yields the exact banded block structure
\begin{align}
 P_kK_YP_k&=a_kP_kM_YP_k,
 &a_k&=\frac{(2k+1)^2-10}{(2k-1)(2k+1)(2k+3)},\label{eq:diagblock}\\
 P_{k+2}K_YP_k&=e_kP_{k+2}M_YP_k,
 &e_k&=\frac3{(2k+1)(2k+3)(2k+5)}.
 \label{eq:offdiagblock}
\end{align}
Here $M_Y$ is multiplication by $Y$; all other blocks vanish, apart from adjoints. The degree-zero diagonal block is zero.

Since $e_0=1/5$, $a_2=1/7$, and $e_2=1/105$, the identities above give
\[
 K_Y\one=\frac15Y,\qquad
 \ip Y{K_YY}=\frac17\int Y^3\dd\mu,\qquad
 P_{\geq3}K_YY=\frac1{105}P_4(Y^2).
\]
The last identity and \eqref{eq:Y4P4} yield $\|P_{\geq3}K_YY\|_2=(6/7)/105=2/245$. These identities give \eqref{eq:block-bound1} and \eqref{eq:block-bound2}.

\subsection*{Diagonal estimates and weighted Schur bound}
The block formula \eqref{eq:diagblock}--\eqref{eq:offdiagblock} leaves only diagonal blocks and couplings between degrees $k$ and $k+2$. We first bound the diagonal blocks through angular momentum operators; a weighted Schur estimate will then control the full operator on the sum of degrees at least three.
Let $J_i=-i(\omega\times\nabla)_i$ denote angular momentum, and set
$J_A=\sum_{i,j}A_{ij}J_iJ_j$. On degree $k$,
\begin{equation}\label{eq:JA}
 P_kM_YP_k=-\frac{2}{(2k-1)(2k+3)}J_A.
\end{equation}
One proof is to use the harmonic decomposition
$(x^\top Ax)h=h_{k+2}+|x|^2h_k+|x|^4h_{k-2}$ for homogeneous harmonic $h$ and compute $\Delta$ and $\Delta^2$ of the product. It gives
\[
 h_{k-2}=\frac{A:\nabla^2h}{(2k+1)(2k-1)},\quad
 h_k=\frac2{2k+3}\left(Ax\cdot\nabla h-\frac{|x|^2}{2k-1}A:\nabla^2h\right),
\]
which is equivalent to \eqref{eq:JA}.

For a complex unit harmonic $h$, the matrix
$M_{ij}=\operatorname{Re}\ip{J_ih}{J_jh}$ is positive semidefinite, with trace $k(k+1)$ and eigenvalues no larger than $k^2$. At this sum, the maximum of $\tr(M^2)$ occurs at eigenvalues $k^2,k,0$, giving
$\tr(M^2)\leq k^4+k^2$. Cauchy--Schwarz for trace-free matrices and \eqref{eq:JA} then yield
\begin{equation}\label{eq:diagestimate}
 \|P_kM_YP_k\|\leq
 \frac{2\sqrt5\,k\sqrt{k^2-k+1}}{(2k-1)(2k+3)}<\frac{\sqrt5}{2}
 \qquad(k\geq1).
\end{equation}
The degree-three block admits the improved estimate
\begin{equation}\label{eq:k3better}
 \|P_3M_YP_3\|\leq\frac{\sqrt5}{3}.
\end{equation}
To see it, diagonalize $A=\operatorname{diag}(a,b,c)$. Degree-three harmonics split into the line generated by $x_1x_2x_3$ and three two-dimensional coordinate-parity spaces. On
$h=x_1(\alpha x_1^2+\beta x_2^2+\gamma x_3^2)$, $3\alpha+\beta+\gamma=0$, the matrix of $J_A$ in $(\beta,\gamma)$ coordinates is
\[
 \begin{pmatrix}-7a-8b&2(c-a)\\2(b-a)&-7a-8c\end{pmatrix},
\]
whose eigenvalues are $-3a\pm3\sqrt{25-4a^2}$. Their absolute values are at most $15\sqrt5/2$; the other parity spaces obey the same estimate, which establishes \eqref{eq:k3better} through \eqref{eq:JA}.

Let $d_k=\|P_kK_YP_k\|$ and $b_k=\|P_{k+2}K_YP_k\|$. From \eqref{eq:diagblock}--\eqref{eq:k3better} and $\|Y\|_\infty\leq\sqrt5$,
\begin{equation}\label{eq:dbounds}
 d_3\leq\frac{13\sqrt5}{315},\qquad
 d_k<\frac4{45}\ (k\geq4),\qquad
 b_k\leq\sqrt5\,e_k\ (k\geq3).
\end{equation}
For the middle bound, substitution into \eqref{eq:diagblock}--\eqref{eq:diagestimate} and exact squaring show $d_4,d_5<4/45$; for $k\geq6$ use
$a_k<(2k+1)^{-1}$ and $\sqrt5/2<9/8$ to obtain $d_k<9/(8(2k+1))\leq9/104<4/45$.

We now apply a weighted Schur bound to $P_{\geq3}K_YP_{\geq3}$, giving degree three weight $w_3=4/3$ and every other degree weight one. On each off-diagonal pair use
$2xy\leq(w_\ell/w_k)x^2+(w_k/w_\ell)y^2$. The degree-three weighted row sum is at most
\[
 \sqrt5\left(\frac{13}{315}+\frac34e_3\right)
 =\frac{617\sqrt5}{13860}=M_*.
\]
We verify the remaining rows explicitly. The diagonal estimate
\eqref{eq:diagestimate}, together with the exact values of $a_4,a_5$,
gives
\[
 d_4\leq\frac{568\sqrt{65}}{53361}<\frac9{100},\qquad
 d_5\leq\frac{370\sqrt{105}}{50193}<\frac2{25}.
\]
Write $b_k\leq\sqrt5\,e_k$ as in \eqref{eq:dbounds}. The degree-four
row is bounded by
\[
 d_4+b_4<\frac9{100}+\frac1{125}=\frac{49}{500},
 \qquad b_4\leq\frac{3\sqrt5}{1287}<\frac1{125}.
\]
The weight $w_3=4/3$ affects only the coupling between degrees three
and five. For the degree-five row,
\[
 d_5+\frac43 b_3+b_5
 <\frac2{25}+\frac7{500}+\frac1{250}
 =\frac{49}{500},
\]
since $\frac43 b_3\leq4\sqrt5/693<7/500$ and
$b_5\leq3\sqrt5/2145<1/250$. Finally, for $k\geq6$ the weights are
all one and the $k$-th row sum is $d_k+b_{k-2}+b_k$. By the bound
$d_k<9/[8(2k+1)]$ established above and the fact that $e_k$
decreases with $k$,
\[
\begin{split}
 d_k+b_{k-2}+b_k
 &<\frac9{104}
  +3\sqrt5\left(\frac1{1287}+\frac1{3315}\right)\\
 &<\frac9{104}
  +\frac{27}{4}\left(\frac1{1287}+\frac1{3315}\right)
 =\frac{9123}{97240}<\frac{49}{500},
\end{split}
\]
where $\sqrt5<9/4$ was used. Thus every row except degree three
is strictly below $49/500$. The degree-three row is bounded by the
previously displayed constant $M_*$, and
\[
 \frac{49}{500}<M_*<\frac1{10}.
\]
Both comparisons follow by squaring. In particular, all the weighted
row sums are bounded by $M_*$. For completeness, this last Schur
estimate can be read directly as
a bound on the quadratic form. Write $G=\sum_{k\geq3}G_k$ with
$G_k=P_kG$; the band structure \eqref{eq:diagblock}--\eqref{eq:offdiagblock}
implies
\[
 |\ip G{K_YG}|
 \leq\sum_{k\geq3}d_k\|G_k\|_2^2
     +2\sum_{k\geq3}b_k\|G_k\|_2\|G_{k+2}\|_2.
\]
For each off-diagonal summand use, with positive weights $w_k$,
\[
 2b_k\|G_k\|_2\|G_{k+2}\|_2
 \leq b_k\frac{w_{k+2}}{w_k}\|G_k\|_2^2
       +b_k\frac{w_k}{w_{k+2}}\|G_{k+2}\|_2^2.
\]
With $w_3=4/3$ and $w_k=1$ for $k\geq4$, the coefficient of each
$\|G_k\|_2^2$ is exactly the weighted row sum estimated above.
It is at most $M_*$, so
\[
 |\ip G{K_YG}|\leq M_*\sum_{k\geq3}\|G_k\|_2^2
 =M_*\|G\|_2^2.
\]
The argument first applies to finite harmonic sums; boundedness of
$K_Y$ and Parseval extend it by density to every admissible $G$. This is exactly \eqref{eq:block-bound3} and completes the proof of Proposition~\ref{prop:blocks}.

\section*{Acknowledgement}
Generative AI tools were used as an aid in the development and presentation of this work during
the preparation of this manuscript. The authors take full responsibility for the mathematical content of this
manuscript, including all statements, proofs, computations, and conclusions.

\end{document}